\documentclass[10pt]{amsart}

\newif\ifprint
\printfalse 

\def\setdrawbox#1#2#3{
\pdfximage{#2} 
\setbox0=\hbox{\pdfrefximage\pdflastximage} 
\drawx=#1\wd0
\ifdim\drawx>\hsize\drawx=\hsize\fi
\pdfximage width \drawx {#2}
\setbox\drawing=\vbox{\offinterlineskip\pdfrefximage\pdflastximage\kern 0pt}
\drawx=\wd\drawing
\drawy=\ht\drawing
\ngap=0pt \sgap=0pt \wgap=0pt \egap=0pt
\setbox0=\vbox{\offinterlineskip \box\drawing \ifgridlines\drawgrid\drawx\drawy\fi #3}%
\setbox\drawing=\vbox{\kern\ngap\hbox{\kern\wgap\box0\kern\egap}\kern\sgap}
}
\def\obsoletedrawbox#1#2#3{\vbox{
 \setbox\drawing=\vbox{\offinterlineskip\epsfbox{#2.eps}\kern 0pt}
 \drawbp=\epsfurx
 \advance\drawbp by-\epsfllx\relax
 \multiply\drawbp by #1
 \divide\drawbp by 100
 \drawx=\drawbp truebp
 \ifdim\drawx>\hsize\drawx=\hsize\fi
 \epsfxsize=\drawx
 \setbox\drawing=\vbox{\offinterlineskip\epsfbox{#2.eps}\kern 0pt}
 \drawx=\wd\drawing
 \drawy=\ht\drawing
 \ngap=0pt \sgap=0pt \wgap=0pt \egap=0pt
 \setbox0=\vbox{\offinterlineskip
 \box\drawing \ifgridlines\drawgrid\drawx\drawy\fi #3}
 \kern\ngap\hbox{\kern\wgap\box0\kern\egap}\kern\sgap}}

\def\drawsmash#1{\relax
\ifmmode%
\typeout{Warning:Macro drawsmash used in math mode}%
\fi%
\setbox0=\hbox{#1}\ht0=0pt\dp0=0pt\box0}

\newcount\figcount
\newbox\drawing
\newcount\drawbp
\newdimen\drawx
\newdimen\drawy
\newdimen\ngap
\newdimen\sgap
\newdimen\wgap
\newdimen\egap

\newif\ifgridlines
\newbox\figtbox
\newbox\figgbox
\newdimen\figtx
\newdimen\figty

\newdimen\bwd
\def\hhline#1{\vbox{\drawsmash{\hbox to #1{\leaders\hrule height \bwd\hfil}}}}
\def\vvline#1{\hbox to 0pt{%
 \hss\vbox to #1{\leaders\vrule width \bwd\vfil}\hss}}

\def\clap#1{\hbox to 0pt{\hss#1\hss}}
\def\vclap#1{\vbox to 0pt{\offinterlineskip\vss#1\vss}}

\def\hstutter#1#2{\hbox{%
 \setbox0=\hbox{#1}%
 \hbox to #2\wd0{\leaders\box0\hfil}}}

\def\vstutter#1#2{\vbox{
 \setbox0=\vbox{\offinterlineskip #1}
 \dp0=0pt
 \vbox to #2\ht0{\leaders\box0\vfil}}}

\def\crosshairs#1#2{
 \dimen1=.002\drawx
 \dimen2=.002\drawy
 \ifdim\dimen1<\dimen2\dimen3\dimen1\else\dimen3\dimen2\fi
 \setbox1=\vclap{\vvline{2\dimen3}}
 \setbox2=\clap{\hhline{2\dimen3}}
 \setbox3=\hstutter{\kern\dimen1\box1}{4}
 \setbox4=\vstutter{\kern\dimen2\box2}{4}
 \setbox1=\vclap{\vvline{4\dimen3}}
 \setbox2=\clap{\hhline{4\dimen3}}
 \setbox5=\clap{\copy1\hstutter{\box3\kern\dimen1\box1}{6}}
 \setbox6=\vclap{\copy2\vstutter{\box4\kern\dimen2\box2}{6}}
 \setbox1=\vbox{\offinterlineskip\box5\box6}
 \drawsmash{\vbox to #2{\hbox to #1{\hss\box1}\vss}}}

\def\boxgrid#1{\rlap{\vbox{\offinterlineskip
 \setbox0=\hhline{\wd#1}
 \setbox1=\vvline{\ht#1}
 \drawsmash{\vbox to \ht#1{\offinterlineskip\copy0\vfil\box0}}
 \drawsmash{\vbox{\hbox to \wd#1{\copy1\hfil\box1}}}}}}

\def\drawgrid#1#2{\vbox{\offinterlineskip
 \dimen0=\drawx
 \dimen1=\drawy
 \divide\dimen0 by 10
 \divide\dimen1 by 10
 \setbox0=\hhline\drawx
 \setbox1=\vvline\drawy
 \drawsmash{\vbox{\offinterlineskip
 \copy0\vstutter{\kern\dimen1\box0}{10}}}
 \drawsmash{\hbox{\copy1\hstutter{\kern\dimen0\box1}{10}}}}}

\def\figtext#1#2#3#4#5{
 \setbox\figtbox=\hbox{#5}
 \dp\figtbox=0pt
 \figtx=-#3\wd\figtbox \figty=-#4\ht\figtbox
 \advance\figtx by #1\drawx \advance\figty by #2\drawy
 \dimen0=\figtx \advance\dimen0 by\wd\figtbox \advance\dimen0 by-\drawx
 \ifdim\dimen0>\egap\global\egap=\dimen0\fi
 \dimen0=\figty \advance\dimen0 by\ht\figtbox \advance\dimen0 by-\drawy
 \ifdim\dimen0>\ngap\global\ngap=\dimen0\fi
 \dimen0=-\figtx
 \ifdim\dimen0>\wgap\global\wgap=\dimen0\fi
 \dimen0=-\figty
 \ifdim\dimen0>\sgap\global\sgap=\dimen0\fi
 \drawsmash{\rlap{\vbox{\offinterlineskip
 \hbox{\hbox to \figtx{}\ifgridlines\boxgrid\figtbox\fi\box\figtbox}
 \vbox to \figty{}
 \ifgridlines\crosshairs{#1\drawx}{#2\drawy}\fi
 \kern 0pt}}}}

\def\setext#1#2#3{\figtext{#1}{#2}10{#3}}

\def\wtext#1#2#3{\figtext{#1}{#2}0{.5}{#3}}
\def\etext#1#2#3{\figtext{#1}{#2}1{.5}{#3}}
\def\ntext#1#2#3{\figtext{#1}{#2}{.5}1{#3}}
\def\stext#1#2#3{\figtext{#1}{#2}{.5}0{#3}}

\def\hpad#1#2#3{\hbox{\kern #1\hbox{#3}\kern #2}}
\def\vpad#1#2#3{\setbox0=\hbox{#3}\dp0=0pt\vbox{\kern #1\box0\kern #2}}

\def\stack#1#2#3{\vbox{\offinterlineskip
 \setbox2=\hbox{#2}
 \setbox3=\hbox{#3}
 \dimen0=\ifdim\wd2>\wd3\wd2\else\wd3\fi
 \hbox to \dimen0{\hss\box2\hss}
 \kern #1
 \hbox to \dimen0{\hss\box3\hss}}}

\def\hexp#1{%
 \setbox0=\hbox{${}^{#1}$}%
 \hbox to .5\wd0{\box0\hss}}

\newlength\gx
\newlength\gy

\def\overgrid#1#2#3#4{\offinterlineskip%
\setlength{\dimen0}{#1}%
\setlength{\dimen1}{#2}%
\divide\dimen0 by #3%
\divide\dimen1 by #4%
\setbox0=\hhline{#1}%
\setbox1=\vvline{#2}%
\rlap{\drawsmash{\vbox{\offinterlineskip%
\copy0\vstutter{\kern\dimen1\box0}{#4}}}}%
\drawsmash{\hbox{\copy1\hstutter{\kern\dimen0\box1}{#3}}}%
}%

\newcommand\labelcapdraw[5]{
\begin{figure}[!ht]
 \refstepcounter{equation}\label{#1}
 \setdrawbox{#3}{#4}{#5}
 \centerline{\ifgridlines\boxgrid\drawing\fi\box\drawing}
 \centerline{{\scshape Figure\nobreakspace\theequation:} \quad {#2}}%
\end{figure}%
}

\usepackage{xspace}
\usepackage{color}
\usepackage[pdftex]{graphicx}
\usepackage{array}

\usepackage{indentfirst,pict2e}
\usepackage{amsfonts}
\usepackage[leqno]{amsmath}
\usepackage{amsthm}
\usepackage{latexsym}

\usepackage{xr-hyper}

\def\docpdftitle{Gr\"{o}bner techniques for low degree Hilbert stability}
\usepackage[
	hyperindex,
	pagebackref,
	pdftex,
	pdftitle={\docpdftitle},
	pdfauthor={Ian Morrison and Dave Swinarski},
	pdfdisplaydoctitle,
	pdfpagemode=UseNone,
	breaklinks=true,
	extension=pdf,
	bookmarks=false,
	plainpages=false,
	colorlinks,
	linkcolor=linkblue,
	citecolor=linkblue,
	urlcolor=linkred,
	pdfmenubar=true,
	pdftoolbar=true,
	pdfpagelabels,
	pdfpagelayout=SinglePage,
	pdfview=Fit,
	pdfstartview=Fit
]{hyperref}

\ifprint
	\definecolor{linkred}{rgb}{0.2,0.2,0.2} 
	\definecolor{linkblue}{rgb}{0.2,0.2,0.2} 
\else
	\definecolor{linkred}{rgb}{0.7,0.2,0.2}
	\definecolor{linkblue}{rgb}{0,0.2,0.6}
\fi

\newcommand{\neturl}[1]{\href{#1}{{\sffamily{\texttt{#1}}}}}
\newcommand{\neturltilde}[2]{\href{#1}{{\sffamily{\texttt{#2}}}}}

\long\def\commentout#1{}
\definecolor{authornote}{rgb}{0,0.6,0} 
\def\authorsnote#1{}

\usepackage[backrefs,msc-links,nobysame]{amsrefs}
\renewcommand{\eprint}[1]{#1}

\makeatletter

\def\bib@div@mark#1{%
 \@mkboth{{#1}}{{#1}}%
	}
\def\print@backrefs#1{%
 \space\SentenceSpace$\leftarrow$\csname br@#1\endcsname
}
\catcode`\'=11 
\renewcommand{\PrintAuthors}[1]{%
 \ifx\previous@primary\current@primary
  \sameauthors\@empty
 \else
  \def\current@bibfield{\bib'author}%
		  \PrintNames{}{}{\scshape #1}%
 \fi
}
\catcode`\'=12
\def\MRhref#1#2{%
 \begingroup
  \parse@MR#1 ()\@empty\@nil%
  \href{\MR@url}{\texttt{\@tempd\vphantom{()}}}%
  \ifx\@tempe\@empty
  \else
   \ \href{\MR@url}{\texttt{(\@tempe)}}%
  \fi
 \endgroup
}%
\def\MR#1{%
 \relax\ifhmode\unskip\spacefactor3000 \space\fi
 \begingroup
  \strip@MRprefix#1\@nil
  \edef\@tempa{\@nx\MRhref{MR\@tempa}{\@tempa}}%
 \@xp\endgroup
 \@tempa
}
\makeatother

\numberwithin{equation}{section}

\newtheoremstyle{claim}
{}
{}
{\itshape}
{}
{\bfseries}
{}
{1em}
{}
\theoremstyle{claim}
\newtheorem{theorem}[equation]{Theorem}
\newtheorem{theorem*}{Theorem}
\newtheorem{definition}[equation]{Definition}

\newtheorem{lemma}[equation]{Lemma}
\newtheorem{corollary}[equation]{Corollary}

\newtheorem{proposition}[equation]{Proposition}
\newtheorem{criterion}[equation]{Criterion}
\newtheorem{pseudoalgorithm}[equation]{Pseudo-Algorithm}

\newtheoremstyle{explain}
{}
{}
{}
{}
{\scshape}
{}
{.5em}
{}
\theoremstyle{explain}
\newtheorem{example}[equation]{Example}
\newtheorem{remark}[equation]{Remark}

\newif\ifarxiv
\arxivtrue
\newcommand{\codesamplelink}[2]{\cite{MorrisonSwinarskiCodeSamples}*{\hyperref{http://www.math.uga.edu/~davids/gs/codesamples.pdf}{codesample}{#1}{{\color{linkred}\texttt{CodeSample~\ifarxiv{#2}\else\ref*{#1}\fi}}}}}

\newcommand{\R}{\mathbb{R}}
\newcommand{\Z}{\mathbb{Z}}

\newcommand{\Q}{\mathbb{Q}}
\newcommand{\G}{\mathbb{G}}
\newcommand{\Gm}{\G_m}
\newcommand{\Grass}{\mathbf{Gr}}
\newcommand{\inI}{\operatorname{in}}

\newcommand{\red}{\mbox{\scriptsize{red}}}
\newcommand{\Pro}{\mathbb{P}}
\newcommand{\Hilb}{\mathbb{H}}

\newcommand{\MUm}{\mbox{\texttt{MUm}}}
\newcommand{\Sym}{\operatorname{Sym}}
\newcommand{\Span}{\operatorname{Span}}

\newcommand{\Aut}{\operatorname{Aut}}
\newcommand{\State}{\operatorname{State}}

\newcommand{\Chow}{\operatorname{Chow}}
\newcommand{\Stab}{\operatorname{Stab}}

\newcommand{\mtwo}{\texttt{Macaulay2} \cite{Macaulay}}
\newcommand{\mtwonocite}{\texttt{Macaulay2}}

\newcommand{\polymakenocite}{\texttt{polymake}}
\newcommand{\gfan}{\texttt{gfan} \cite{Gfan}}
\newcommand{\gfannocite}{\texttt{gfan}}
\newcommand{\Magma}{\texttt{MAGMA}~\cite{Magma}}
\newcommand{\Magmanocite}{\texttt{MAGMA\xspace}}
\newcommand{\GAP}{\texttt{GAP}~\cite{Gap}}
\newcommand{\StatePolytope}{\texttt{StatePolytope}\xspace}

\newcommand{\Phat}{{\ensuremath{\widehat{P}}}}
\newcommand{\Qhat}{{\ensuremath{\widehat{Q}}}}
\newcommand{\Rhat}{{\ensuremath{\widehat{R}}}}
\newcommand{\Ihat}{{\ensuremath{\widehat{I}}}}
\newcommand{\muhat}{{\ensuremath{\widehat{\mu}}}}
\newcommand{\Hilbhat}{{\ensuremath{\widehat{\Hilb}}}}

\newcommand{\SL}{\operatorname{SL}}
\newcommand{\GL}{\operatorname{GL}}
\newcommand{\PGL}{\operatorname{PGL}}
\newcommand{\diag}{\operatorname{diag}}
\newcommand{\Hom}{\operatorname{Hom}}
\newcommand{\M}{\ensuremath{\mathcal M}}

\newcommand{\Mbar}{\ensuremath{\overline{\M}}}
\newcommand{\Mibar}[1]{\ensuremath{\Mbar_{#1}}}

\newcommand{\Mgbar}{\ensuremath{\Mibar{g}}}

\newcommand{\m}{\ensuremath{M}}

\newcommand{\mbar}{\ensuremath{\overline{\m}}}
\newcommand{\mibar}[1]{\ensuremath{\mbar_{#1}}}

\newcommand{\mgbar}{\ensuremath{\mibar{g}}}

\newcommand{\HO}[2]{H^0\bigl(#1,#2\bigr)}
\newcommand{\hO}[2]{h^0\bigl(#1,#2\bigr)}
\newcommand{\bary}{{\mathbf{0}}}
\renewcommand{\o}{{\mathcal O}}

\newcommand{\poly}{{\mathcal P}}
\newcommand{\thst}[2]{\ensuremath{{#1}^{\mathrm{#2}}}}
\newcommand{\lbpow}[2]{#1^{\otimes #2}}

\newcommand{\resp}[1]{{\upshape [resp: }{#1}{\upshape{\kern0.1em]}}\xspace}
\newcommand{\wiman}[1]{{\mathcal W}_{#1}}
\newcommand{\oneps}{$1$--ps\xspace}

\newcommand{\restrictedto}[2]{%
\ensuremath{{#1\hskip 1pt{\vrule height 7.2pt depth 3.6pt}\hskip0.75pt}%
\raisebox{-1.8pt}[0pt][0pt]{\ensuremath{\null_{#2}}}}}

\begin{document}

\def\normallinespread{1.13}
\def\normalspread{\linespread{\normallinespread}\normalfont\selectfont}
\def\abslinespread{0.8}
\def\absspread{\linespread{\abslinespread}\footnotesize\selectfont}
\lineskiplimit=-0.08cm

\pagenumbering{arabic}

\title{Gr\"{o}bner techniques for low degree Hilbert stability}

\author{Ian Morrison}
\address{Department of Mathematics\\ Fordham University\\ Bronx, NY 10458}
\email{morrison@fordham.edu}

\author{David Swinarski}
\address{Department of Mathematics\\ University of Georgia\\ Athens, GA 30602}
\email{davids@math.uga.edu}


\abstract{{\absspread We give a method for verifying, by a symbolic calculation, the stability or semistability with respect to a linearization of fixed, possibly small, degree $m$, of the Hilbert point of a scheme $X \in {\mathbb P}(V)$ having a suitably large automorphism group. We also implement our method and apply it to analyze the stability of bicanonical models of certain curves. Our examples are very special, but they arise naturally in the log minimal model program for} $\overline{\mathcal M}_g$}{\absspread. In some examples, this connection provides a check of our computations; in others, the computations confirm predictions about conjectural stages of the program.}}


\maketitle

\normalspread 

\section{Introduction} \label{introsection}

We analyse the Hilbert stability of bicanonical models of certain curves $X$ of small genus with suitably large automorphism groups with respect to linearizations of fixed small degree $m$. Our examples are very special, but they have geometrically interesting applications discussed below.

Our analysis has two main novelties. First, we give a method for deducing the stability, always with respect to $\SL(V)$, of the Hilbert point of a subscheme $X$ of $\Pro(V)$, from a symbolic calculation of certain state polytopes. Even the possibility of such a reduction for Hilbert points of subschemes of large codimension is new. The key hypothesis we use is that $X$, as a subscheme of $\Pro(V)$, is \emph{multiplicity free} (\ref{multfreedef}): the multiplicity, in the natural representation of $\Aut(X)$ on $V$, of every irreducible representation is either $0$ or $1$. In our examples, $V = \HO{X}{\lbpow{\omega_X}{2}}$ and we say that $X$ is \emph{bicanonically multiplicity free} (\ref{Xbicanonicallymultiplicityfree}). Most of our examples are certain special hyperelliptic curves $\wiman{g}$, called Wiman curves, that are well known in the literature on curves with automorphisms \cite{BreuerCharacters}, and nodal curves that are joins of two or more Wiman curves. 

The second novelty of our examples is that they allow us to handle fixed values of the linearization degree $m$. Indeed, the values of $m$ arising in our applications are not merely fixed but quite small, typically $6$ or less. In contrast, existing approaches, such as those pioneered by Gieseker in \cites{GiesekerGlobal, GiesekerTata, GiesekerCIME} have an asymptotic character and verify Hilbert stability only with respect to linearizations of sufficiently large degree $m$.

The bicanonical curves $X$ and small degrees $m$ in our main examples are chosen because quotients of loci in the bicanonical Hilbert schemes $\Hilb$ in question are predicted to yield new log minimal models of the moduli spaces of stable curves. For further details on this connection, see 7.5 of \cite{MorrisonGIT} and the references cited there. A disclaimer is in order here. Examples in which $X$ is smooth show the non-emptiness of stable loci of interest in the log minimal model program in small genus but are far from producing the desired quotients. This paper deals only with our methods for checking stability and with our examples; it discusses the construction of such quotients only in outlining our plans for further work.  

In addition to providing applications of our examples, the log minimal model program makes, by indirect arguments, very specific predictions about the degrees in which the bicanonical Hilbert points of certain reducible and singular $X$ will be stable, strictly semistable and unstable. A number of our examples involve curves of these types. For some of these, we use completed stages of the minimal model program as a check on our calculations. For others, our calculations verify the program's predictions exactly, providing further evidence for them.

Our approach combines the hypothesis of multiplicity freeness with theorems of Kempf on worst destabilizing \oneps's to reduce checking stability for the full group $\SL(V)$ to checking stability with respect to a distinguished maximal torus $T$ (Corollary~\ref{slhilbertstabilityfrompolytope}).  There is an easy naive algorithm for checking this symbolically, but its complexity makes it impractical except in very simplest cases. By adapting results of Bayer and the first author on state polytopes, we give, in Corollary~\ref{BMcorollary}, an algorithm efficient enough that we are able to handle examples arising in our intended applications. 

Working with small degree $m$ is a sword that cuts both ways. On the one hand, the $m$ we work with are well below the bounds that ensure various standard uniformity hypotheses for ideals of points of $\Hilb$, even those that are deformations of smooth subschemes. A typical example is that the degree $m$ graded pieces of the homogenous ideals do not yield the embedding of $\Hilb$ as a closed subscheme of a Grassmannian needed to linearize the $\PGL(V)$-action. We address these complications by replacing $\Hilb$ with a multigraded Hilbert scheme $\Hilbhat$ in the sense of Haiman and Sturmfels \cite{HaimanSturmfels}.

On the other hand, our algorithm is only practical when computing state polytopes in fairly low degrees. It involves computing \emph{all} the monomial initial ideals $X$ (in the coordinates giving the special torus $T$) and requires a Gr\"obner basis calculation for \emph{each} initial ideal. In fact, as the genus of $X$---and hence the bicanonical embedding dimension---increased, we were often unable to carry even these low degree calculations to completion because there are simply too many such ideals. To understand such examples, we use several additional, somewhat ad-hoc tricks. 

The first involves a Monte Carlo strategy that computes a random sub-polytope of the state polytope by computing some random initial ideals. If $X$ is Hilbert stable and we are fortunate, this sub-polytope provides a proof of stability. This approach can never prove that $X$ is unstable, but we are able to do this, when necessary, by educated guesswork. Geometry---in our examples, analogies with completed stages of the log minimal model program---often suggests what a destabilizing $\lambda$ should be, and such a guess can be checked by a \emph{single} Gr\"obner basis computation. Finally, the Parabola Trick (Proposition~\ref{hhlemma}) uses ideas of Hassett, Hyeon and Lee~\cite{HassettHyeonLee} to deduce stability of Hilbert points of smooth curves in low degrees not accessible to our calculations from their stability in even lower degrees.

Here is the plan of the rest of the paper. The details of our multigraded setup for Hilbert points is given in  Section~\ref{hilbertpointsection} and of the results on state polytopes we need are extended to this setting in Section~\ref{statepolysection}. Section~\ref{Kempfsection} reviews Kempf's results on worst one-parameter subgroups and explains how, for multiplicity free $X$, they reduce checking stability to calculations with state polytopes. The Monte Carlo version and the Parabola Trick are outlined in section~\ref{additionaltrickssection}. Section~\ref{wimansection} recalls facts about Wiman curves and pluricanonical equations of hyperelliptic curves needed to set up the \mtwonocite\ calculations for our examples and section~\ref{predictionssection} summarizes stability properties that are either known from or predicted by the log minimal model program for $\mgbar$. Our examples and how they fit with the log minimal model program are reviewed in Section~\ref{resultssection}. Finally, we close by listing some ideas for future work. 

\subsection*{Accessing raw and commented source code}
Our calculations are carried out in \mtwo\ using the \StatePolytope package of the second author that calls the packages \gfan\ and \polymakenocite\ \cites{GawrilowJoswigPolymake, Polymake} to compute intermediate results. Some preliminary calculations are performed in \Magma\ or \GAP. The source code of our routines and detailed output from many calculations are posted on the second author's webpage at \neturltilde{http://www.math.uga.edu/~davids/gs/gs.html}{http://www.math.uga.edu/$\sim$davids/gs/gs.html}. For the convenience of the reader, we have collected annotated code snippets in a separate paper \cite{MorrisonSwinarskiCodeSamples}, posted in the same directory, to which \href{http://www.math.uga.edu/~davids/gs/codesamples.pdf}{{\color{linkred}\texttt{Code Sample}}} (hyper)references here point.

\subsection*{Acknowledgements}
It is a great pleasure to thank Brendan Hassett, David Hyeon, and Yongnam Lee, whose work got us interested in these GIT questions, and Dave Bayer, Johan de Jong, Bill Graham, Anders Jensen, Julius Ross, Greg Smith, David Smyth, Mike Stillman and Bernd Sturmfels for helpful discussions. The first author wishes to acknowledge support from a Fordham University Faculty Fellowship during the early stages of this project. The second thanks Sonja Mapes for her expert instruction on programming in \mtwonocite, and Dan Grayson for suggesting many improvements to the \StatePolytope package.  

\section{Hilbert points and state polytopes} \label{hilbertpointsection}

\subsection*{Parameter schemes adapted to low degrees} 
Fix an $(N+1)$-dimensional vector space $V$ over an algebraically closed ground field $K$, and a set of coordinates $\{x_{0},\ldots,x_{N}\}$ identifying $V$ with $K^{N+1}$ and the homogeneous coordinate ring of $\Pro(V)$ with $S := K[x_{0},\ldots,x_{N}]$. Fix also a Hilbert polynomial $P$ of degree $r$ and let $\Hilb$ be the Hilbert scheme of $r$-dimensional subschemes $X\subset \Pro(V)$ with Hilbert polynomial $P$. 

The goal of this section is to define state polytopes of such subschemes $X$---or of their homogeneous ideals $I \subset S$---and to recall their connection to the stability of the Hilbert point of $X$ with respect to the action of $\SL(V)$ induced by the natural action on $\Pro(V)$. Both of these notions depend on the choice of the degree $m$ that is used to linearize this action.  To make uniform sense of either the Hilbert point or the state polytope for \emph{all} $X$ having a fixed Hilbert polynomial $P$---that is, over the whole of $\Hilb$---it is necessary to take $m$ larger than the Gotzmann number (\cite{Gotzmann}) $m_P$ for $P$.  An easy calculation using the formulae there gives these numbers for Hilbert polynomials of curves:  
\begin{lemma}  \label{gotzmannformula}
The Gotzmann number of $P(t)= at+b$ is $m_{P} = \binom{a}{2} + b$. 
\end{lemma}

However, the applications we have in mind to stability problems arising in the log minimal model program for $\mgbar$ (cf. \cites{HassettHyeonLogCanonical, HassettHyeonFlip}) require us to work with a fixed degree $m < m_P$. The main goal of this section is to outline how to transfer the standard constructions to this setting.  This is most conveniently achieved by using the multigraded Hilbert schemes constructed by Haiman and Sturmfels \cite{HaimanSturmfels}. In doing this we have treated general $r$, since doing so entails no additional complications, but for the applications cited above, we will specialize to the case $r=1$ of curves.

We begin with a definition of convenience.
\begin{definition} \label{lnice}
An $r$-dimensional subscheme $X$ of $\Pro(V)$ with ideal sheaf $I$ is called \emph{$\ell$-nice} if
\begin{enumerate}
\item The natural map $V^\vee \to \Gamma(X, \o_X(1))$ is an isomorphism. 
\item $\o_X$ is $(\ell-1)$-regular.
\item $I_X$ is $\ell$-regular.
\end{enumerate}
\end{definition}

This list of properties comes from a similar list of hypotheses for certain statements in  \cite{HassettHyeonFlip}.  We note one minor change: all of the subschemes we start with are pure $r$-dimensional (though their specializations under various \oneps may not be) and so we have omitted this condition in Definition \ref{lnice}.

The first hypothesis may be viewed more geometrically as saying that $X$ is embedded in $\Pro(V)$ by a complete non-degenerate linear series. The second hypothesis implies that for $m \ge \ell$, the sheaf $\o_X(m)$ has no higher cohomology, and hence that its Hilbert polynomial $P(m)$ computes $\hO{X}{\o_X(m)}$. Likewise, the third hypothesis implies that the restriction maps $S_{m} \to \HO{X}{\o_X(m)}$ are surjective for all $m \ge \ell$ and that $I_X$ is generated be elements of degree at most $\ell$.

We denote by $\Hilb_{\ell}$ the $\ell$-nice locus in the Hilbert scheme $\Hilb$ of subschemes of $\Pro(V)$ with Hilbert polynomial $P$. Fix an $\ell$-nice subscheme $X$. We let $\Rhat(m) = \dim_K(S_m) = \binom{m+N}{m}$ and $\Qhat_{\ell}(m)= \dim_K(I_m)$ for $m\ge \ell$ and $\Qhat_{\ell}(m)= 0$ for $m<\ell$. In other words, $\Qhat_{\ell}$ is the Hilbert function of the ideal $\Ihat_{\ell}$ given by truncating $I$ in degrees below $\ell$. As usual, we can recover $I$ from any $\Ihat_{\ell}$ by saturating. Our hypotheses imply that $\Ihat_{\ell}$ is generated in degree exactly $\ell$. Finally, let $\Phat_{\ell}(m) = \Rhat(m)-\Qhat_{\ell}(m)$. This is a truncation of the Hilbert function of $X$ and only equals $P(m)$ for $m \ge m_P$. 

We denote by $\Hilbhat_{\ell}$ the multigraded Hilbert scheme of ideals in $S$ with Hilbert function $\Phat_{\ell}$ and denote by $[I]$ the point of $\Hilbhat_{\ell}$ determined by the ideal $I$. By \cite{HaimanSturmfels}*{Corollary 1.2}, $\Hilbhat$ is a projective scheme representing the functor of \emph{locally free} families of such ideals and hence is equipped with a universal family. Their Lemma 4.1 identifies $\Hilbhat_{m_P}$ with the usual Hilbert scheme $\Hilb$ and, if $\ell < m_P$, then truncation up to degree $m_P$ gives a map $i_{\ell}: \Hilbhat_{\ell} \to \Hilb$. 

A few cautions are in order here. First, the $\ell$-nice locus in $\Hilbhat_{\ell}$ is only locally closed, and it need not even be dense---there may be entire components of $\Hilbhat_{\ell}$ containing no $\ell$-nice ideals. 

Second, while $i_{\ell}( \Hilbhat_{\ell})$ is closed in $\Hilb$ and $i_{\ell}$ is injective on the $\ell$-nice locus, the map $i_{\ell}$ need not be an embedding. This pathology has its origin in the  fact that the ideals parameterized by $\Hilbhat_{\ell}$ need not be saturated, even in degrees above $\ell$ where they are not truncated. For example, if $\Hilb$ contains a point $X'$ whose (saturated) ideal $I'$ satisfies  $\dim_K(I'_{\ell}) > \dim_K(I_{\ell})$ and $\dim_K(I'_m) = \dim_K(I_m)$ all $m > \ell$, then every choice of a $\dim_K(I_{\ell})$-dimensional subspace of $I'_{\ell}$  determines an ideal $I'' \in \Hilbhat_{\ell}$ mapping to $X'$. Such examples can be found, for example, with $\Hilb$ the Hilbert scheme of twisted cubics (cf. Example~\ref{twistednoninjectiveg}). The upshot is that we cannot replace the ideal $I$ parameterized by a point of $\Hilbhat_{\ell}$ by the subscheme $X$ it determines unless we know that the degree $\ell$ truncation of the saturation of $I$ has Hilbert function exactly $\Phat_{\ell}$, as we do, by definition, over the $\ell$-nice locus.

\subsection*{The Hilbert matrix} 

Henceforth we fix values of $\ell$ and $m \ge \ell$.  In our applications, we often take $\ell=2$.   We begin with two remarks designed to lighten our notation. First, since $m \ge \ell$, $\Phat_{\ell}(m)$ depends only on $m$, so we can and will omit the subscript $\ell$s used above. Second, we introduce many objects depending on our choice of $m$ in this section, but when there is no risk of confusion, we will omit the $m$ to simplify notation in later sections.

Let $W_m = \bigwedge^{\Phat(m)} S_m$ and let $\Grass_m \subset \Pro(W)$ be the  Pl\"ucker embedding of the  Grassmannian $\Grass_m := \Grass\bigl(\Phat(m), \Rhat(m)\bigr)$ of $\Phat(m)$-dimensional quotient spaces of $S_m$. There is a Pl\"ucker map $g_m: \Hilbhat_{\ell} \to \Grass_m $ sending $[I]$ to $S_m/I_m$. The map $g_{m}$ has closed image, but need not be injective: for example, $g_2$ has the same value on the monomial ideals \ref{third} and \ref{fourth} in Example~\ref{twistednoninjectiveg} and on the ideals \ref{seventh} and \ref{eighth}. If, however, $I$ is generated in degrees at most $m$---in particular, for points in the $l$-nice locus---$g_{m}([I])$ does determine $I$.

We want to describe homogeneous coordinates $y_A$ of $g_m([I]) \in \Pro(W)$ in a form usable in tools like \mtwonocite. This is most conveniently and concretely done by working with the \emph{subspace} $I_m $ of $S_m$ rather than the quotient $S_m/I_m$, and using it to define $m$-Hilbert matrices $M_{I,m}$. First let $\mathcal{B}_{m} = \{x^j\}$ be the monomial basis of $S_m$ with a fixed ordering. Then let $\mathcal{C}_{m}(I) =\{p_i, i = 1, \ldots , \Phat(m) \}$ be any ordered basis of $I_m$ and let $M_{I,m}$ be the $\Phat(m) \times \Rhat(m)$ matrix whose $\thst{ij}{th}$ entry is the coefficient of the monomial $x^j$ in the equation $p_i$.  The Pl\"ucker \emph{coordinates} $y_A$ of $I_m$ are then simply the  $\Qhat(m) \times \Qhat(m)$ minors of $M_{I,m}$---one for each Pl\"ucker \emph{set} $A$ of $\Qhat(m)$ of the monomials $\mathcal{B}_{m}$.  As in the discussion on page 211 of \cite{BayerMorrison}, if $M'_{I,m}$ is the matrix associated to a second basis $\mathcal{C}'_{m}(I)$ and $E$ is the associated change of basis matrix, then $M' = EM$ and, for all $A$, $y'_A= \det(E) y_A$. Hence,
\begin{enumerate}
	\item The point $g_m([I])$ of $\Pro(W)$ defined by the collection of $y_A$ is independent of the choice of $\mathcal{C}_{m}(I)$. 
	\item Whether or not any individual $y_A$ vanishes at $g_m([I])$ is likewise independent of this choice. 
	\item We may always make this choice so that $M_{I,m}$ is in echelon form. 
\end{enumerate}

\begin{example} 
	For a monomial ideal, we may take the basis $\mathcal{C}_{m}$ to be monomial, too, and then the Hilbert matrix is particularly simple: it will have exactly one $1$ in each row and be $0$ otherwise.  Thus, for a given $m$, there is exactly one nonzero Pl\"{u}cker coordinate, given by the Pl\"{u}cker set $A=\mathcal{C}_{m}$.
\end{example}

\begin{example} \label{pluckerexample}
	Consider the ideal $I$ of two distinct points in $\Pro^{2}$. For instance $P=(1,2,3)$ and $Q=(5,1,-4)$.  Let $a,b,c$ be the coordinates on $\Pro^{2}$.  Then we can view  $I$ as $(c-3a,b-2a) * (a-5b,c+4b)$ and take
\begin{displaymath}
\mathcal{C}_{2} = [2a^2-11ab+5b^2,8ab-4b^2+2ac, 3a^2-15ab-ac+5bc, 12ab+3ac-4bc-c^2]\,.
\end{displaymath}
Ordering $\mathcal{B}_{S_{2}}$ as $[a^{2},ab,ac,b^2,bc,c^2]$, we get:
\begin{equation}
M_{I, 2} = \left(
\begin{array}{cccccc}
2 & -11 & 0 & 5 & 0 & 0 \\
0 & 8 & 2 & -4 & -1 & 0 \\
3 & -15 & -1 & 0 & 5 & 0 \\
0 & 12 & 3 & 0 & -4 & -1 
\end{array}
\right)
\end{equation}
Then the Pl\"{u}cker point of $M_{I,2}$ is given by the following point, in which we have indexed the Pl\"ucker sets by the pair of monomials \emph{omitted} to save space:

\medskip
{
\setlength{\tabcolsep}{2pt}
\begin{tabular}{m{.5cm}m{.08cm}m{.5cm}m{.08cm}m{.5cm}m{.08cm}m{.7cm}m{.08cm}m{.5cm}m{.08cm}m{.5cm}m{.08cm}m{.5cm}m{.08cm}m{.5cm}m{.08cm}m{.5cm}m{.08cm}m{.5cm}m{.08cm}m{.5cm}m{.08cm}m{.7cm}m{.08cm}m{.5cm}m{.08cm}m{.5cm}m{.08cm}m{.9cm}}
$\widehat{12}$ && $\widehat{13}$ && $\widehat{14}$ && $~\widehat{15}$ && $\widehat{16}$ && $\widehat{23}$ && $~\widehat{24}$ && $\widehat{25}$ && $~\widehat{26}$ && $~\widehat{34}$ && $\widehat{35}$ && $~\widehat{36}$ && $\widehat{45}$ && $\widehat{46}$ && $~\widehat{56}$ \\
45 &:&-95 &:&99 &:& -154 &:&209 &:& 55 &:& -18 &:& 38 &:& -13 &:& -83 &:& 108 &:& -228 &:& 22 &:& 55 &:& -132~.
\end{tabular}
}

\vskip6pt
\noindent Alternatively, the Pl\"ucker coordinates can be computed in \mtwonocite\ \codesamplelink{pluckercode}{2}. 
\end{example}

\section{Stability and state polytopes} \label{statepolysection}

\subsection*{$T$-states and $T$-state polytopes} \label{mudef}
We next want to focus on the action of $\SL(V) \cong \SL(N+1)$ on $W:= \bigwedge^{P(m)} \Sym^{m} V$.  The Hilbert--Mumford criterion says that $w \in W$ is $\SL(V)$ stable if and only if $w$ is $\lambda$-stable for every $1$-parameter subgroup (henceforth \oneps) $\lambda: \Gm \rightarrow \SL(V)$.   If, in terms of a basis of $V$  with respect to which $\lambda$ diagonalizes as $\diag(t^{r_0}, \cdots , t^{r_N})$, the point $w$ has coordinates $(w_0, \ldots, w_N)$, we set
\[
\mu^{L}(w,\lambda) :=  -\min \{ r_{i}  \mid i \mbox{ such that } w_{i}^{*} \neq 0 \},
\]
and $w$ is $\lambda$-stable if and only if $\mu^{L}(w,\lambda) <0$.  

\begin{remark} \label{musignremark}
A word about the minus sign in the definition of $\mu$. Our preferred sign convention for the index $\mu$ of a Hilbert point $w$ is that of \cite{GiesekerTata},  \cite{BaldwinSwinarski} and \cite{MorrisonGIT} in which we consider the Grassmannian as parameterizing $P(m)$-dimensional quotients of $S_m$, given by restriction to $\HO{X}{\o_X(m)}$, and $w$ is stable if any $\lambda$ acts with negative weight on some non-zero coordinate of $w$.  
	
	The minus sign has been inserted to compensate for the fact that here we will be calculating weights of the action of $\lambda$ on the degree $m$ piece of the ideal $I$ of $X$ which is of dimension $\Qhat(m)$. This, of course, gives rise to a quotient of dimension $\Phat(m)$ and the complement of each Pl\"ucker set $A$ of monomials gives a basis of this quotient. But when we take $m$ to be small, we can no longer identify the quotient with $\HO{X}{\o_X(m)}$ except on the $l$-nice locus, and it therefore seemed easier to us to simply work with $I_m$. This choice has no effect on the notion of $\SL(V)$-stability because the possibility of replacing $\lambda$ by its inverse means that $w$ is stable if and only if we can always find non-zero coordinates of $w$ on which $\lambda$ acts with weights of opposite signs. 
\end{remark}

\begin{remark} \label{mumremark}
When we are considering stability of a Hilbert point $[X]$, we will write $\mu([X],\lambda)(m)$ for the index with respect to the degree $m$ linearization. Replacing the line bundle $L_m$ on $\Hilb$ or $\Hilbhat$ that is being linearized by the degree $m$, as we have implicitly been doing to this point, is harmless because $m$ determines $L_m$ up to powers. 
	
	For fixed choices of $x$ and $\lambda$, the index is represented by a polynomial in $m$ when $m$ is sufficiently large---for example, when $m \ge l$ if $X$ and its specialization under $\lambda$ are both $l$-nice (cf. \cite{HassettHyeonFlip}*{\S 3.7, 3.8}, and Section~\ref{parabolatricksubsection} below). In our examples, we will often be interested in the roots of this polynomial and in other aspects of the dependence of the index $\mu$ on the degree $m$, so it is convenient to use a notation that highlights this dependence. 
	
	In some examples arising out of the log minimal model program for $\mgbar$, we will want to consider rational values of $m$. We can interpret what our calculations say in such cases without going into the somewhat complicated process of constructing a linearization with respect to a rational value; for details on these, see \cite{DolgachevHu} or \cite{ThaddeusVGIT}. 
	
	All we will use is that the weights of points with respect to such a linearization interpolate those with respect to integral linearizations in the following sense: if $[X]$ and its $\lambda$ specialization are $l$-nice, then $\mu([X],\lambda)(m)$ is computed by a polynomial in $m$ for any rational $m \ge l$. This allows us first to determine this polynomial from values at integral $m$, and then to evaluate it to find $\mu([X],\lambda)(m)$ at the rational values of $m$ that are of interest.
\end{remark}

The connection with Gr\"obner theory comes via another way of expressing the stability of $w$ with respect to the maximal torus $T$ of $\SL(V)$ determined by a choice of basis $B$ of $V$.  Any character $\chi\in \Hom(T, \Gm)$ of $T$ may be written
\[ \chi \bigl( \diag(d_{0},\ldots,d_{N}) \bigr) = \prod_{i=0}^{N} d_{i}^{z_{i}}.
\]
where the $z_i$ are integers, determined, since we are in $\SL(V)$, up to a common shift. Further, any representation $W$ of $T$  decomposes into a direct sum of character eigenspaces $W_{\chi}$, where $w\in W_{\chi}$ if and only if $t\cdot w = \chi(t) w$ for all $t \in T$.   

Define the $T$-\textit{state} $\State_T(w)$ of $w$ to be the set of characters for which the eigencomponent $w_{\chi}$ of $w$ is non-zero, and define the $T$-\textit{state polytope} $\poly_T(w)$ to be the convex hull of $\State_T(W)$ in $\Hom(T, \Gm)$\footnote{This is the state polytope of \cite{BayerMorrison}*{\S 2} defined entirely in the fixed degree $m$, as opposed to that of Sturmfels \cite{SturmfelsGrobner}*{Theorem 2.5} which is the Minkowski sum of the former for all degrees up to $m$.}.

The group of $1$-parameter subgroups of $T$ is dual to its character group: $\lambda \cdot \chi$ is the $\lambda$-weight of $\chi$---the power of $t$ determined by the homomorphism $\chi\circ\lambda:\Gm \to \Gm$.  Viewing $\lambda$ as giving a linear functional on $\Hom(T, \Gm)$, we may rephrase our discussion of the Numerical Criterion as saying that $w \in W$ is stable with respect to a \oneps $\lambda$ in $T$ if and only $w$ has two non-zero eigencomponents $w_{\chi}$ whose character vectors lie on opposite sides of the hyperplane on which $\lambda$ vanishes.  Thus we arrive at the following characterization of GIT stability: 

\begin{criterion} \label{barycentercritgeneral}
A vector $w \in W$ is $T$-stable iff the trivial character lies in the interior of the state polytope and is $T$-strictly semistable iff the trivial character lies on the boundary of the state polytope.
\end{criterion}

To interpret Criterion~\ref{barycentercritgeneral} for Hilbert points, first observe that each eigenspace $(S_m)_{\chi}$ of $S_m$ is spanned by a single $B$-monomial $M$ and, if we normalize the choice of the $z_i$ above by requiring that they sum to $m$, then we may identify the character $\chi$ and the exponent vector of $M$. The Pl\"ucker coordinates $y_A$ on $W$ likewise give an eigenbasis, although the eigenspaces are not necessarily $1$-dimensional. If we now normalize so that the $z_{i}$ sum to $\Qhat(m)m$, then we can identify the corresponding character $\chi_A$ with the sum of the exponent vectors of the $\Qhat(m)$ monomials determined by $y_A$.  For example, in $\bigwedge^{2} \Sym^{2} K[a,b,c,d]$, the wedge product $a^2 \wedge bc$ lies in the weight space corresponding to $W_{(2,1,1,0)}$.   

Monomials and Pl\"ucker coordinates also diagonalize the actions of a \oneps $\lambda$ of $T$ on $S_m$ and $W$. The weight $w_\lambda(M)$ of a monomial $M$ is the sum of the weights of its coordinate factors and the weight $w_\lambda(y_A)$ of a Pl\"ucker coordinate $y_A$ is the sum of the weights of the monomials in it. Moreover, these weights agree with the $\lambda$-weights of the corresponding characters.

Thus, we think of the characters as lying on the hyperplane 
\[ Z_m := \{ z \in \Z^{N+1} \mid \sum_{i=0}^{N} z_{i} = mP(m) \}
\]
This identifies the trivial character with the point in $\Q^{N+1}$ having all coordinates equal to  $\frac{m\Qhat(m)}{N+1}$. In the sequel, we will denote this point by $\bary_m$ and call it the \textit{barycenter} of $Z_m$. 

To simplify two notations that we will use frequently, we write $\State_{T,m}(I)$ and $\poly_{T,m}(I)$ for the $T$-state and the $T$-state polytope of $g_m([I])$, omitting the $T$ when no confusion is possible.  

\begin{criterion} \label{barycentercrithilbert}
The \thst{m}{th}-Hilbert point $g_m([I])$ of an ideal $I$ is $T$-stable iff $\bary_m$ lies in the interior of $\poly_{T,m}(I)$ and is $T$-strictly semistable iff $\bary_m$ lies on the boundary of $\poly_{T,m}(I)$.
\end{criterion}

Note that Criteria \ref{barycentercritgeneral} and \ref{barycentercrithilbert} only test $T$-stability.  In Corollary~\ref{slhilbertstabilityfrompolytope}, we will identify conditions under which we can extend this to $\SL(V)$-stability. 

\begin{example}
	The $\thst{m}{th}$-state of a monomial ideal $I$ is a single point, since there is only one nonzero Pl\"{u}cker coordinate. Unless this point equals $\bary_m$, the \thst{m}{th}-Hilbert point of $I$ is unstable.
\end{example}

\begin{example}
	If $X$ a hypersurface of degree $d$ in $\Pro^N$, we may take $d=m$---so $P(m)=1$---and suppress the exterior power in $W$. Both the characters appearing in the decomposition of $W$ and its Pl\"ucker coordinates are then indexed by monomials $\prod_{i=0}^N x_i^{z_i}$ of degree $d$, and, viewed as lying in the plane $\sum_{i=0}^N z_I=d$, form the $\thst{d}{th}$ subdivision of an $N$-simplex.

\labelcapdraw{pcpdfl}{Degree $3$ state polytope for a cuspidal plane cubic with ideal $I= \langle x^2z=y^3\rangle$}{0.6}{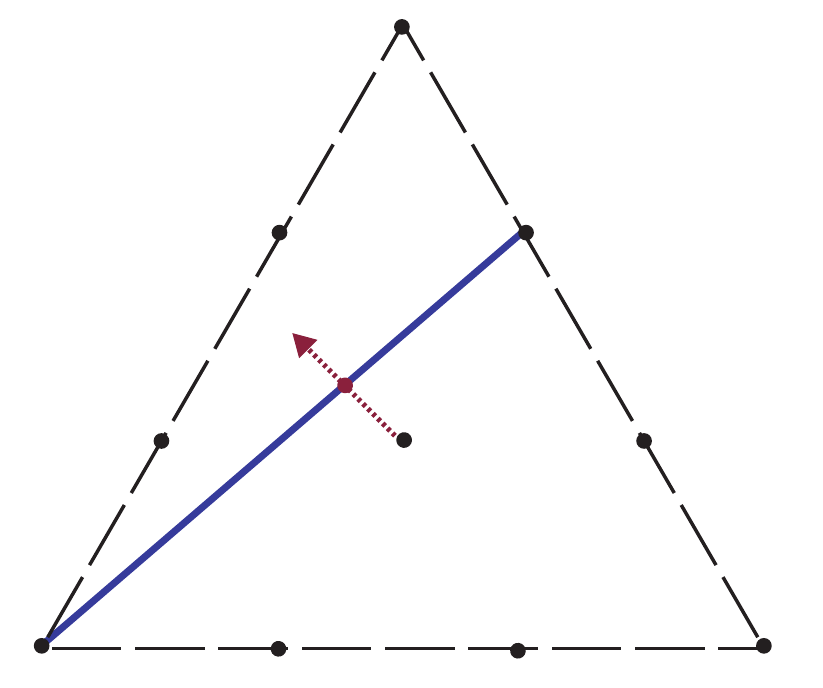}{
		\etext{.05}{.0}{{$y^3$}}
		\stext{.50}{1.01}{$x^3$}
		\wtext{.94}{.0}{$z^3$}
		\wtext{.68}{.68}{{$x^2z$}}
		\setext{.37}{.54}{{$\lambda$}}
		\ntext{.49}{.33}{$\bary$}
		\wtext{.53}{.31}{{\tiny$=(1,1,1)$}}
		\wtext{.44}{.455}{{{\tiny$\leftarrow\kern-2pt \frac{1}{8}(10,9,5)$}}}
		\wtext{.6}{.85}{{{$\lambda=\diag(t^4,t,t^{-5})$}}}
		\wtext{.20}{.20}{{$\leftarrow \poly_{T,3}(I)$}}
		}%
Figure~\ref{pcpdfl} shows this situation for a cuspidal plane cubic $X$ with equation $x^2z=y^3$ which is unstable with respect to the \oneps $\lambda$ shown. The set of characters appearing in the decomposition is are indicated by dots and the simplex that is their convex hull is the outlined triangle. The state polytope is the line segment joining the two monomials with non-zero coefficients in the equation.
	
	These monomials both have weight $3$ with respect to the \oneps $\lambda$ given by $\lambda(t)=\diag(t^4,t,t^{-5})$ and hence this Hilbert point is unstable. The instability is reflected in the fact that $\poly_{T,3}([X])$ does not contain~$\bary$. 
	
	A generic hypersurface in $\Pro^{2}$ would have a two-dimensional state polytope.  The degeneracy of $\poly_{T,3}([X])$ reflects the fact that this cuspidal cubic has a $\Gm$ action (however, it is not normal). But adding an $x^3$ term to the equation, making the state polytope the upper sub-triangle subtended by $\poly_{T,3}([X])$ would not affect the instability.	
\end{example}

\begin{example} We return to Example~\ref{pluckerexample} and the ideal of $[1:2:3] \cup [5:1:-4] \subset \Pro^{2}$.  Every character with non-zero eigenspace contains a nonzero Pl\"{u}cker coordinate, and the state polytope is the two-dimensional hexagon pictured in Figure~\ref{twopointfig}. Here the barycenter (indicated by the central solid circle) has coordinates $(\frac{8}{3}, \frac{8}{3}, \frac{8}{3})$. 
\begin{center}
\setlength{\unitlength}{.35cm}
\begin{picture}(5,6.)(-3.5,-3.5)
\multiput(-2,1.41)(2,0){3}{\circle{.2}}
\multiput(-3,0)(2,0){4}{\circle{.2}}
\multiput(-2,-1.41)(2,0){3}{\circle{.2}}
\multiput(-1,-2.83)(2,0){2}{\circle{.2}}
\put(-2,1.41){\line(1,0){4}}
\put(-3,0){\line(100,141){1}}
\put(-3,0){\line(100,-141){2}}
\put(3,0){\line(-100,141){1}}
\put(1,-2.83){\line(100,141){2}}
\put(-1,-2.83){\line(1,0){2}}
\put(-2.25,1.61){\tiny 431}
\put(-0.25,1.61){\tiny 422}
\put(1.75,1.61){\tiny 413}
\put(-3.25,.2){\tiny 341}
\put(-1.25,.2){\tiny 332}
\put(-.2,-.6){\tiny $\bullet$}
\put(.75,.2){\tiny 323}
\put(2.75,.2){\tiny 314}
\put(-2.25,-1.21){\tiny 242}
\put(-0.25,-1.21){\tiny 233}
\put(1.75,-1.21){\tiny 224}
\put(-1.25,-2.63){\tiny 143}
\put(.75,-2.63){\tiny 134}
\end{picture}
\refstepcounter{equation}\label{twopointfig}
\centerline{{\scshape Figure\nobreakspace\theequation:} \quad $\State_2(I)$ for two general points in $\Pro^{2}$}%
\end{center}
\end{example}

\subsection*{Vertices of state polytopes and initial ideals} 

The number of Pl\"ucker coordinates grows quickly as the number of variables, the number of generators of the ideal, and $m$ grow.  Thus it is impractical to compute the state polytope from definitions for all but the very simplest examples.  The following results, modeled closely on analogous statements in \cite{BayerMorrison}, allow us to handle larger examples by giving a procedure for finding the vertices of $\poly_{T,m}(I)$ that avoids the need to deal with interior Pl\"ucker coordinates.

Any \oneps $\lambda$ in $T$ yields a partial order $\ge_\lambda$ on monomials:
\begin{displaymath}
	M \ge_\lambda M' \iff w_{\lambda}(M) \ge w_{\lambda}(M')\,.
\end{displaymath}
Since the weights $w_i$ of $\lambda$ on $V$ are integers there will always be ties in large degree. But in any fixed degree $m$, $\ge_\lambda$ will give a total order for all $\lambda$ not lying on a finite collection of hyperplanes. If so, we say that $\lambda$ is $m$-generic. We will say that $\lambda$ is generic if it is $m$-generic for $\ell \le m \le m_P$. 

\begin{lemma} \label{BMlemma} For any $I$ in $\Hilbhat$ and any generic \oneps $\lambda$, there is a unique Pl\"ucker set $A_{\lambda}$ of $\Qhat(m)$ monomials such that:
	\begin{enumerate}
		\item $y_{A_{\lambda}}$ is non-zero at $g_m([I])$.
		\item If $y_{A'}$ is any other Pl\"ucker coordinate non-zero at $g_m([I])$, then $w_{\lambda}(y_A) > w_{\lambda}(y_{A'})$.
	\end{enumerate} 
	Moreover, if $M_{I,m}$ is a $m$-Hilbert matrix for $I$ in echelon form, then the monomials in $A_{\lambda}$ span the $>_\lambda$-initial ideal $\inI_{>_\lambda}(I)$ of $I$ in degree $m$.
\end{lemma}
\begin{proof} This is the content of Lemma 3.3 and Corollary 3.4.(ii)  of \cite{BayerMorrison} and the proofs given there apply verbatim in our situation.
\end{proof}

\begin{definition} For any generic \oneps $\lambda$, we let $\chi_{\lambda}= \chi_{A_{\lambda}}$. In other words, $\chi_{\lambda}$  is the character given by summing the exponent vectors of the $\Qhat(m)$ monomials in $\inI_{>_\lambda}(I)_m$. By (1) of Lemma~\ref{BMlemma}, this character is an element of $\State_m(I)$.
\end{definition}

\begin{theorem} \label{BMtheorem} For any $m$-generic \oneps $\lambda$, the character $\chi_{\lambda}$ is a vertex of the state polytope $\poly_{T,m}(I)$. Conversely, if $\chi$ is any vertex of $\poly_{T,m}(I)$, then the eigenspace $W_{\chi}$ is one dimensional and is spanned by the Pl\"ucker coordinate  $y_{A_{\lambda}}$ for some $m$-generic $\lambda$. In particular, $\chi = \chi_{\lambda}$.
\end{theorem}
\begin{proof}
The inequality in (2) of Lemma~\ref{BMlemma} shows that $\sum_{i=0}^{N} w_i z_{i} = w_{\lambda}(y_A)$ is a supporting hyperplane ($\chi_{\lambda}$ lies on it and all other $\chi'$ in $\State_m(I)$ lie on the negative side of it) and hence proves the first claim.  Conversely, any supporting hyperplane $\sum_{i=0}^{N} w_i z_{i} = b$ to $\chi$ may be perturbed so that the coefficients of its normal are the set of weights $w_i$ of a generic \oneps $\lambda$.  But then any Pl\"ucker coordinate $y_A$ lying in the $\chi$-eigenspace satisfies the conditions defining $y_{A_{\lambda}}$ in Lemma~\ref{BMlemma}. The Lemma therefore implies that there is a unique such Pl\"ucker coordinate and that $\chi = \chi_{\lambda}$. The second claim follows.
\end{proof}

We note that, in general, the dimension of $W_{\chi}$ will be quite large. Already in Figure~\ref{twopointfig}, the three interior characters have $2$-dimensional eigenspaces. 

Theorem~\ref{BMtheorem} is a weaker version of Theorem 3.1 of \cite{BayerMorrison} which shows that if $m \ge m_P$, then the set of vertices of $\poly_{T,m}(I)$ is canonically bijective to the set of initial ideals of $I$. For the small degrees that we are treating here where the map $g_m$ from an ideal to its degree $m$ graded piece is not injective, a surjection from initial ideals to vertices is all that we can hope for---and all we need for our applications. 

\begin{example}\label{twistednoninjectiveg} 
Let $X$ be the twisted cubic in $\Pro^3$ with ideal $I= \langle ac-b^2,ad-bc,bd-c^2\rangle $.  Then $X$ has eight initial ideals (see \codesamplelink{usingstatepolytopepackage}{3}):
\begin{enumerate}
	\item $\langle  bd, ad, ac\rangle $ 
	\item $\langle  c^2, ad, ac\rangle $ 
	\item \label{third} $\langle  c^2, bc, ac, a^2d\rangle $ 
	\item \label{fourth} $\langle  c^2, bc, b^3, ac\rangle $ 
	\item $\langle  c^2, bc, b^2\rangle $ 
	\item $\langle  bd, b^2, ad\rangle $ 
	\item \label{seventh} $\langle  bd, bc, b^2, ad^2\rangle $ 
	\item \label{eighth}$\langle  c^3, bd, bc, b^2\rangle $
\end{enumerate}
$\State_2(I)$ has six vertices: initial ideals \ref{third} and \ref{fourth} agree in
degree $2$, as do initial ideals \ref{seventh} and \ref{eighth}.  For any $m\geq 3$,  $\State_m(I)$ has eight vertices.	
\end{example}

By \cite{BayerThesis}*{Proposition~1.8}, given any multiplicative total order $>$ we can find a \oneps $\lambda$ such that $>$ and  $>_\lambda$ agree up to degree $m$. Hence,

\begin{corollary} \label{BMcorollary} The state polytope $\poly_{T,m}(I)$ is the convex hull of the set of $\chi_A$ as $A$ runs over all Pl\"ucker sets that are bases for the degree $m$ graded piece of some initial ideal of $I$.
\end{corollary}

Conveniently, Anders Jensen's program \gfan\ computes the set of initial ideals of $I$.  Thus, if we compute the $\thst{m}{th}$ state of each initial ideal for sufficiently large $m$, we will have the state polytope.  This is what the \mtwo\ package \texttt{StatePolytope} does. 

We will not use the following geometric characterization of $A_{\lambda}$ but have found it helpful in thinking about the preceding results. The action of $\SL(V)$ on $V$ induces actions on the homogeneous polynomials of each degree on $V$ and hence on Hilbert scheme $\Hilbhat$ and on the Grassmannian $\Grass_m$ for which the map $g_m$ is equivariant. Since $\Hilbhat$ is projective, we can define an ideal $J$ giving a point of $\Hilbhat$ by
\begin{displaymath}
	[J] := \lim_{t\to 0} \lambda(t) \cdot [I]\,.
\end{displaymath}
Lemma~\ref{BMlemma} says that $y_{A_{\lambda}}$ is the unique Pl\"ucker coordinate that is non-zero at $g_m([J])$ and hence that  $g_m([J]) = g_M([\inI_{>_\lambda}(I)])$. But all these arguments apply equally to any other degree between $\ell$ and $m_P$ so that, in all these degrees, $J$ and $\inI_{>_\lambda}(I)$ are equal. Hence, $J=\inI_{>_\lambda}(I)$ in degrees above $\ell$.

\section{Kempf's theory of the worst \oneps} \label{Kempfsection}

Let $w$ be a point of an $\SL(V)$ representation $W$.  Already on page 64 in the first edition of \cite{GIT}, Mumford conjectured that if $w$ is unstable, there is a worst destabilizing \oneps $\lambda$ as measured by the index $\mu(w, \lambda)$. In this section, we review the proof of this conjecture by Kempf~\cite{KempfInstability} and Rousseau~\cite{Rousseau} but, to simplify, deal only with the linear situation we need in our applications. We follow the treatment of Kempf, which contains some complementary results that allow us to reduce the $\SL(V)$-stability of points $w \in W$ with suitably large stabilizer $\Stab(w)$ to their $T$-stability for a special torus determined by this stabilizer.

We begin by reviewing some of the background of Kempf's arguments. First, some easy covariance properties. 
\begin{lemma} \label{muinvariances} For any $g\in \SL(V)$,
\begin{enumerate}
	\item $\State_T(g\cdot w) = g\State_T(w)g^{-1}$. 
	\item \label{muconjugation} $\mu(w,\lambda) = \mu(g\cdot w, g\cdot \lambda\cdot g^{-1})$. 
\end{enumerate}
\end{lemma}
\begin{proof}
	The first statement is Lemma~3.2.d) of \cite{KempfInstability} and the second follows immediately from it.
\end{proof} 

Next, let $F_{\lambda}$ be the $\lambda$ weight filtration on $V$ and $P_{\lambda}$ be the parabolic subgroup of block upper triangular matrices in $\SL(V)$ preserving $F_{\lambda}$. Equivalently, $P_{\lambda}$ consists of those $p \in \SL(V)$ for which the limit $\lim_{t \to 0}\lambda(t) \cdot p \cdot \lambda^{-1}(t)$ exists. 

\begin{lemma} \label{Pinvariances} If $g\in \SL(V)$ and $p \in P_{\lambda}$ then, \begin{enumerate}
	\item \label{Pconjugation} $P_{g\lambda g^{-1}} = g P_{\lambda} g^{-1}$.
	\item $p \in P_{\lambda} \iff P_{p\lambda p^{-1}} = P_{\lambda}$.
	\item \label{mupcovariance} If $p \in P_{\lambda}$, $\mu(w,\lambda) = \mu( w, p \lambda p^{-1}) = \mu(p^{-1}\cdot w,\lambda)$. Hence, $\mu(w,\lambda) = \mu(p\cdot w,\lambda)$.
\end{enumerate}
\end{lemma} 
\begin{proof}
	The first statement follows directly from the characterization of $P_{\lambda}$ in terms of limits and then the second follows because any parabolic subgroup is its own normalizer.
	
	The last assertion is trickier. Our proof follows that of Lemma~3.3.e) from~\cite{KempfInstability}. The final equality follows from the first two by inverting $p$. By Lemma~\ref{muinvariances}(\ref{muconjugation}), the second equality follows from the first. For this, the key point is the following claim: if, as $t \to 0$, $\lambda(t) p^{-1} \lambda^{-1}(t) \to p_0^{-1}$ and $t^{-r}\lambda(t)\cdot w \to w_0$, then $t^{-r} p\lambda(t)p^{-1}\cdot w \to p p_0^{-1}\cdot w_0$. This follows because
\begin{displaymath}
	t^{-r}p\lambda(t)p^{-1}\cdot w = p\,\bigl(\lambda(t)p^{-1} \lambda^{-1}(t)\bigr)\,\bigl(t^{-r}\lambda(t)\cdot w \bigr)\to p\,\, p_0^{-1}\cdot w_0
\end{displaymath}
Since $\mu(w,\lambda)$ is the largest $r$ such that $\lim_{t\to 0}t^{-r}\lambda(t)\cdot w$ exists, the claim shows that if $p \in P_{\lambda}$, then $\mu(w,\lambda) \le \mu(w,p\lambda p^{-1})$ and, by symmetry, (\ref{mupcovariance}) follows. 
\end{proof} 

Replacing $\lambda(t)$ by $\lambda_k(t) := \lambda(t^k)$ for any positive integral $k$ scales all weights by $k$ without affecting their signs. Thus stability with respect to $\lambda$ and $\lambda_k$ are equivalent, but $\mu(w, \lambda_k) = k \mu(w, \lambda)$. We want to normalize the Hilbert-Mumford index $\mu$ index to obtain a measure of ``badness'' that agrees on $\lambda$ and $\lambda_k$. To do this, choose a conjugation-invariant norm $||\cdot||$ on one parameter subgroups---for $\SL(V)$, we can take 
$||\lambda|| := \left(\sum_{i=0}^N w_i^2\right)^{\kern-1pt\lower3pt\hbox{\scriptsize$\frac{1}{2}$}}$---
and define $\muhat(w,\lambda) := \frac{\mu(w,\lambda)}{||\lambda||}$. We also define $\overline{\mu}(w) := \sup_{\lambda} \muhat(w,\lambda)$. A priori, it is not clear either that $\overline{\mu}(w)$ is finite or, if it is, that this sup is achieved. A \emph{worst} $\lambda$ for $w$ is one for which the function $\muhat(w,\lambda)$ achieves this maximum value.

In the case when $X$ is a representation $W$, Theorem 3.4 of \cite{KempfInstability} says that:

\begin{theorem} If $w $ is an unstable point of $W$, then
\begin{enumerate}
	\item There is an indivisible \oneps $\lambda$ such that, if $\lambda'$ any other \oneps, then $\muhat(w,\lambda) \ge \muhat(w,\lambda')$. Hence $\overline{\mu}(w)$ is finite and equal to $\muhat(w,\lambda)$.
	\item \label{worstunicity} The indivisible $\lambda'$ for which $\muhat(w,\lambda') =\muhat(w,\lambda)$ are exactly those for which $\lambda' = p^{-1} \lambda_w p$ for some $p \in P_{\lambda}$. In particular, $P_{\lambda'} = P_{\lambda}$ and we can write $P_w$ for $P_{\lambda}$.
	\item The set of all $\lambda'$ as in~(\ref{worstunicity}) is a principal homogeneous space under the unipotent radical of $P_{w}$, and every maximal torus $T$ of $P_{w}$ contains a unique such $\lambda'$. 
\end{enumerate} 
\end{theorem}

In view of Lemma~\ref{muinvariances}(\ref{muconjugation}) and Lemma~\ref{Pinvariances}(\ref{mupcovariance}), we can informally summarize this result as saying that worst one-parameter subgroups exist and are as unique as possible. The complementary result we need is:

\begin{proposition}[\cite{KempfInstability}*{Corollary 3.5}] \label{Kempfprop} Let $ w \in W$ be an unstable point with associated parabolic subgroup $P_w$. Then $P_{w}$ contains $\Stab_{w}(\SL(V))$.
\end{proposition}
\begin{proof}
	For any $g \in \SL(V)$, $g\cdot w$ is also unstable and hence determines a parabolic subgroup $P_{g\cdot w}$. The Proposition will follow, if we show that $g P_w g^{-1} = P_{g\cdot w}$, because any parabolic subgroup is its own normalizer. 
	
	By Lemma~\ref{muinvariances}(\ref{muconjugation}) and the conjugation invariance of the norm $||\cdot||$, the (indivisible) worst one-parameter subgroups for $g\cdot w$ are exactly the $g$-conjugates of those for $w$. Let $\lambda$ be one of the latter. This gives the middle equality in  $P_{g\cdot w} = P_{g \lambda g^{-1}} = g P_{\lambda} g^{-1} = g P_w g^{-1}$ and the first and last equalities follow from Lemma~\ref{Pinvariances}(\ref{Pconjugation}).
\end{proof}

Kempf applies these results to conclude stability of Chow and Hilbert points of abelian varieties and homogeneous spaces (\cite{KempfInstability}*{Cor. 5.2 and 5.3)}: the representations of the automorphism groups of these varieties are irreducible, so the stabilizer is not contained in any nontrivial parabolic, and hence these must be GIT stable.  

There are very few examples of pluricanonically embedded smooth curves with an automorphism group acting via an irreducible representation. For instance, a full list of canonical curves with this property is found in \cite{BreuerCharacters}*{App. B}; the highest genus example is $g=14$.  Examples are even rarer as $\nu$ increases.  So, Kempf's strategy must be modified if it is to be applied to Hilbert points of curves. 

Here is how we weaken the irreducibility hypothesis.
\begin{definition}\label{multfreedef}
	Fix $w \in W$. We say that $w$ is \emph{multiplicity free} with respect to a finite subgroup $G$ of $\Stab_{\SL(V)}(w)$ if, in the representation of $G$ on $V$, no $G$-irreducible $R$ has multiplicity greater than $1$. When, as in our applications here, $G=\Stab_{\SL(V)}(w)$ we will simply say that $w$ is \emph{multiplicity free}.
\end{definition}

The key consequence of this property is that $V$---indeed, any $G$-invariant subspace $U$ of $V$---has a \emph{canonical} decomposition as a direct sum of $G$-irreducible subrepresentations of $V$. Such a decomposition, of course, exists for any finite $G$ by complete reducibility. But when $w$ is multiplicity free, there is, for each $R$ appearing in $V$, a canonical subrepresentation $U_R$ isomorphic to $R$. Every $U$ is then the direct sum of those $U_R$ for which $R$ occurs in $U$.

\begin{definition} \label{determinesstability}
	We say that a basis $B$ of $V$ or the associated torus $T=T_B$ \emph{determines stability} for $w$ if: 
\begin{enumerate}
	\item There is a subgroup $G \subset \Stab_{\SL(V)}(w)$ such that $w$ is multiplicity free with respect to $G$.
	\item The basis $B$ is the (disjoint) union of bases $B_R$ for each of the $G$-irreducible representations $U_R$ occurring in $V$.
\end{enumerate}
\end{definition}

The justification for this terminology is:

\begin{proposition}\label{multiplicityfreetorus}
	If $T$ determines stability for $w$ and $w$ is $T$-semistable, then $w$ is $\SL(V)$-semistable.
\end{proposition}
\begin{proof}
	We prove that if $w$ is $\SL(V)$-unstable, then $w$ is $T$-unstable, by showing that then $T$ is a torus of $P_w$.
	
	So suppose that $w$ is multiplicity free and unstable. Then Proposition~\ref{Kempfprop} says that $G$ lies in $P_w$ and hence fixes the associated filtration $F$. We can thus write $F$ as a strictly nested sequence $V=U_0 \supset U_1 \supset \cdots \supset U_h \subset \{0\}$ of $G$-invariant subspaces of $V$. Each of these is a direct sum of a subset of the $G$-irreducibles occurring in $V$. Therefore, the basis $B$ is compatible with the filtration $F$ and, in turn $T$ is a torus of $P_w$.
\end{proof}

We now apply this to Hilbert points. Let $X$ be an $l$-nice subscheme of $\Pro(V)$ with ideal $I$ and let $\Aut_V(X) \subset \Aut(X)$ be the subgroup of consisting of elements that act linearly on $V$ fixing $X$. Suppose that $\Aut_V(X)$ is finite and the representation of $\Aut_V(X)$ on $V$ is multiplicity free---in our applications $\Aut_V(X) = \Aut(X)$. For any $m \ge l$, the group $\Aut_V(X)$ lies in the $\SL(V)$-stabilizer of $g_m([I])$, so the pair $\bigl(g_m([I]), \Aut_V(X)\bigr)$ is multiplicity free in the sense of Definition~\ref{multfreedef}, independently of $m$. 

\begin{definition} \label{multfreeX}
	Under these hypotheses of the preceding paragraph, we say that \emph{$X$ is multiplicity free} and that any torus $T$ constructed as in Proposition~\ref{multiplicityfreetorus} \emph{determines stability for $X$}. 
\end{definition}

Combining the Proposition with Criterion~\ref{barycentercrithilbert} gives the first assertion below. The second, which allows us to read off the worst \oneps from the state polytope follows by elementary arguments as in the proof of \cite{KempfInstability}*{Lemma 2.3}.

\begin{corollary}\label{slhilbertstabilityfrompolytope}
\begin{enumerate}
	\item If $T$ determines stability for $X$, then the \thst{m}{th}-Hilbert point $g_m([I])$ of $X$ is $\SL(V)$-stable \resp{$\SL(V)$-strictly semistable} if and only if the barycenter $\bary_m$ lies in the interior \resp{the boundary} of the state polytope $\poly_{T,m}(I)$.
	\item Let $p$ be the proximum to $\bary$ in $\poly_{T}(w)$. If $w$ is $T$-unstable, so $p \not = 0$, then $p-\bary_m$ spans a rational ray and any $T$-worst \oneps has weights lying on this ray.
\end{enumerate}
\end{corollary}

\begin{remark}
Of course, multiplicity free Hilbert points are extremely special. Consider, for example, smooth curves of genus $g \ge 2$. These usually have trivial automorphism group, and so, for any embedding, trivial stabilizer. But even special curves can only be multiplicity free for low degree embeddings as discussed further below. So our strategy can only prove directly the stability of low degree models of special curves, like those arising in our applications here. In practice, the complexity of the computations required blows up very rapidly---for an indication of just how rapidly, see Table~\ref{wimanfourtable}---making it practical to handle even such cases only for small $g$.

On the other hand, by the openness of GIT stability and the coarseness of the Zariski topology, proving that a single smoothable subscheme in any component of the Hilbert scheme is stable proves that a general smooth subscheme on that component is stable. As Gieseker's construction of $\mgbar$, and many others modeled on it (cf.~\cite{MorrisonGIT}), show, such a statement is often enough, when there is a main component containing smooth equidimensional subschemes to allow the construction of a GIT quotient to be completed by indirect arguments. 
\end{remark}

	For the rest of the paper we specialize to the case where $X$ is a curve, though many arguments will continue to apply more generally. To make this switch clear we write $C$ for $X$, continuing to denote its ideal by $I$. In looking for examples in this case, the next step is therefore clear. Find special models $C \subset \Pro(V)$ of curves that are multiplicity free and decide when their  \thst{m}{th}-Hilbert points are stable, at least for small $m$, by computing $\poly_{T,m}(I)$ for some $T$ that determines stability. A natural set of models to consider are pluricanonical ones, since for these $\Aut_V(C) = \Aut(C)$. 
	
\begin{definition}\label{Xbicanonicallymultiplicityfree}
We say that a nodal curve $C$ is \emph{$\nu$-multiplicity free} if its $\nu$-canonical model is multiplicity free. We will mainly be interested in the case  $\nu=2$ when we say that $C$ is \emph{bicanonically multiplicity free}.
\end{definition}

In the sequel, we focus on bicanonical models for two reasons. First, they provide a source of tractable examples. Bicanonical embedding dimensions are small enough both so that multiplicity free examples exist in all genera and so that it is practical to compute the relevant state polytopes when $g$ is sufficiently small. Second, as noted in the introduction, the conjectural next stages in the log  minimal model program of Hassett and Hyeon~\cites{HassettHyeonLogCanonical, HassettHyeonFlip} depend on understanding the stability of bicanonical Hilbert points of degree at or below $6$. 

Finally, MacLachlan \cite{Maclachlan}*{Theorem 4} shows that, if a curve of genus $g$ has an \emph{abelian} automorphism group, its order can be at most $4g+4$. Hence no such curve can be $\nu$-multiplicity free for any $\nu >2$ unless is has very small genus---all the irreducibles have dimension $1$, and the $\nu$-canonical series has dimension $(2\nu-1)(g-1)$. Multiplicity free curves with non-abelian automorphism groups seem likewise to be extremely rare. Since the sum of the squares of the dimensions of the irreducibles equals the order of the group, the sum of the dimensions themselves is typically much less, and this more than compensates for the largely theoretical extra headroom given by Hurwitz' bound of $84(g-1)$ for their orders. Our examples, reviewed in Section~\ref{wimansection}, use a family $\wiman{g}$ of hyperelliptic curves, one in each genus $g$, called Wiman curves whose automorphism groups are \emph{cyclic} of order---surprisingly, in view of Maclachlan's bound---$4g+2$.

\section{Complementary approaches to checking stability} \label{additionaltrickssection}

We are almost ready to describe our applications of the preceding results to check, by symbolic calculations, GIT stability of certain bicanonically multiplicity free curves $C$ of small genus with respect to small $m$ linearizations. Given such an $C$ with ideal $I$, we want to compute the state polytope $\poly_{T,m}(I)$ with respect to a distinguished torus of Corollary~\ref{slhilbertstabilityfrompolytope} and check whether the barycenter $\bary_m$ lies in its interior, boundary or exterior.  But in practice, already in some cases with $g=4$, this plan is impossible to carry out: there are too many initial ideals and we are unable to completely compute $\poly_{T,m}(I)$. In this section, we explain some additional ideas that we use to settle such cases in our examples. 

\subsection*{A Monte Carlo Pseudo-Algorithm} \label{montecarlosubsection}

Our first observation is that we can often check GIT stability without computing the entire state polytope. If we can compute \emph{any} set $\mathcal I$ of initial ideals corresponding, in degree $m$, to a set of characters $\Xi_m \subset \State_m(I)$ such that the convex hull $\overline{\Xi}_m$ contains $\bary_m$, then we know $C$ is $m$-Hilbert semistable---even stable if $\bary_m$ lies in the interior of $\overline{\Xi}_m$. We say such an $\mathcal I$ checks $m$-semistability of~$C$.

\begin{pseudoalgorithm}\label{montecarlopseudo}
	To verify that $C$ is $m$-Hilbert semistable, begin with $\mathcal I = \emptyset$.
	\begin{enumerate} 
		\item \label{monteone}Generate a pseudo-random weight vector $\lambda$.
		\item Add the ideal $\inI_{>_\lambda}(I)$ to $\mathcal I$ and let $\Xi_m$ be the associated set of characters in degree $m$.
		\item If $\overline{\Xi}_m$ contains $\bary_m$, stop. Otherwise, return to Step~(\ref{monteone}).
		\end{enumerate} 
\end{pseudoalgorithm}

This is, of course, only a pseudoalgorithm because it will never terminate if $C$ is actually $m$-Hilbert \emph{un}stable. In fact, even if $C$ is $m$-Hilbert semistable, we cannot be sure it will terminate: we may simply not have generated vectors in enough directions to produce a $\Xi$ that checks semistability. In practice, however, we have not encountered either problem. Guided by the predictions of the log minimal model program, we have been able to apply it only to testing $C$ that actually were $m$-Hilbert semistable. In all the examples we have run, Pseudo-Algorithm~\ref{montecarlopseudo} has produced a $\Xi$ that checks stability quickly, typically in a tiny fraction of the time required for our calculation of the full state polytope to complete---or fail to complete due to hardware and software limitations.

We have also been able to check the $m$-Hilbert instability of $C$ in all cases we expect it by exploiting a basic asymmetry of stability calculations: verifying semistability is hard, but proving instability is easy, if we can guess a destabilizing \oneps. To check such a guess, it suffices to compute an the $\inI_{>_\lambda}$ initial ideal (adding a tie-breaking procedure if necessary). The ideal $I$ is $m$-Hilbert unstable if and only if the monomial basis of this ideal in degree $m$ has negative weight: by Lemma~\ref{BMlemma}, this weight equals $\mu(g_m(I), \lambda)$. The \mtwonocite\ function \texttt{MUm} computes $\mu$ by this method.

In our examples, the log minimal model program suggests what \oneps to try.  It gives a geometric description of those curves whose instability it predicts that suggests both the filtration and the weights of a candidate destabilizing \oneps $\lambda$. These candidates have proven to be destabilizing in all our examples. In practice, we can often even check instability with respect to these $\lambda$ deductively, because $\lambda$-weights have a geometric interpretation. This pattern is familiar to those who have computational experience, symbolic or deductive, with Hilbert stability.

\subsection*{The Parabola Trick} \label{parabolatricksubsection}

For a fixed curve $C\in \Pro(V)$, the complexity of computing the state polytope $\poly_{T,m}(I)$ grows quite rapidly with $m$ because $\binom{m+N}{m}$, the number of monomials of degree $m$, grows like $N^m$ for $N \gg m$. This often means that we can check $\ell$-Hilbert stability for some $\ell$ but not $m$-Hilbert stability for $m$ a somewhat larger, but still small, degree of greater geometric interest. Typically $\ell=2$ and, with applications to the log minimal model program for $\mgbar$ in mind, $m\le 6$. The following proposition sometimes lets us deduce what we want to know from what we can compute.

The key technical tool is a Lemma due to Hassett and Hyeon:

\begin{proposition}[The Parabola Trick] \label{hhlemma} Let $\lambda$ be a non-trivial \oneps.  Suppose that $C$ and its $\lambda$ specialization are both $\ell$-nice.  Then the weight function $\mu([C],\lambda)(m)$ is computed by a quadratic polynomial for $m \ge \ell$, and this polynomial has the form  $a(m-1)(m-r)$ for some rational $a$ and $r$. In particular, any two values of  $\mu([C],\lambda)(m')$ with $m' \ge \ell$ determine $\mu([C],\lambda)(m)$ for all $m\ge \ell$.
\end{proposition}
\begin{proof} For $\ell=2$, this follows from \cite{HassettHyeonFlip}*{Proposition~3.17}, but the same proof works, mutatis mutandi, for any $\ell \ge 2$.
\end{proof}

We apply this trick to check both instability and stability. A typical example of the latter use is:

\begin{corollary}\label{parabolcor} Suppose $C$ is $2$-nice, and that we can find a set  $\mathcal I$ of initial ideals that are generated in degree $\ell$ and such that for all $m$ between $\ell$ and $m_{P}$, the corresponding set $\{ \State_{m}(I_j) \mid I_j \in \mathcal{I} \}$ contains $\bary_m$ in its convex hull. Then $C$ is $m$-Hilbert stable for all $m \geq \ell$.  
\end{corollary}

The Gotzmann number $m_{P}$ grows very quickly with the genus of the curve, and 
therefore it is natural to wonder if perhaps a smaller set of degrees $m$ might 
suffice in the corollary.   The following example dashes any such hopes.

\begin{example} Let $\wiman{4}$ be the genus 4 Wiman curve discussed in Example \ref{wimanfourcalc} below.  Let $I$ be the ideal of this curve under its bicanonical embedding (see \ref{wimanfourexample} below).  The state polytope of $I$ is a subset of $\R^{9}$.  There is a set $\{ I_j \}_{j=1}^{9}$ of nine initial ideals of $I$ that are generated in degrees $\leq 7$.  Hence $I$ and all the $I_j$ are 7-nice.  For $4\leq m \leq 36$, but \emph{not} for $37 \leq m \leq 64$, the set $\{ \State_{m}(I_j) \}_{j=1}^{9}$ contains $\bary_m$ in its convex hull.  Thus we cannot apply the Parabola Trick to conclude stability from this set of initial ideals (although, as we explain in \ref{wimanfourcalc} below, we can find another set of initial ideals that does establish stability).
\end{example}

\begin{remark}\label{chowremark}
	In Sections~\ref{hilbertpointsection} and~\ref{statepolysection}, we emphasized the changes that must be made to the theory as it appears in the literature in order to accommodate low degrees $m$.  But our approach also applies (in any dimension $r$) to degrees above the Gotzmann number $m_P$, and via these may be used to determine Chow stability. The idea is straightforward. Fix a saturated ideal $I$ corresponding to a $r$-dimensional subscheme $X$. Given any \oneps $\lambda$, let $\xi_m$ be the character associated to $\inI_{>_\lambda}(I)$ in degree $m$. For sufficiently large $m$,  $\xi_m$ is represented by a vector of polynomials in $m$ and we may  define a scaled limit $ \xi_{\infty} = \lim_{m \rightarrow \infty } \frac{k!}{m^k} \xi_m$. By \cite{KapranovSturmfelsZelevinsky}*{Theorem 3.3}, the Chow polytope (in \cite{KapranovSturmfelsZelevinsky}, simply the state polytope) $\Chow(I)$ may then be defined either as the convex hull of these scaled limits for all $\lambda$, or as the scaled limit of the degree $m$ state polytopes of $I$. Thus, we may compute Chow polytopes by computing state polytopes in $r+1$ sufficiently large degrees $m$, interpolating in each coordinate, and computing the scaled limit.  
\end{remark}

\section{Background on Wiman Curves} \label{wimansection}

Looking ahead to the next section, our main source of computational examples will be a sequence $\wiman{g}$ of hyperelliptic curves called Wiman curves. In this section, we develop the theoretical background on these curves that we will need in these applications. To begin with, we need to have suitable equations for their pluricanonical models and this part of the story depends only on their being hyperelliptic.

\subsection*{Pluricanonical equations of hyperelliptic curves}

So fix a smooth hyperelliptic curve $C$ and a canonical multiple $\nu$, and let $\phi_{\nu}: C \rightarrow \Pro(V_{\nu})= \Pro^{N_{\nu}}$ be the $\nu$-canonical embedding in $V_{\nu} = \HO{C}{\lbpow{\omega_C}{\nu}}^{\vee}$.  By Riemann-Roch, $N_{\nu}+1 = (2 \nu -1)(g-1)$ and $\phi_{\nu}(C)$ has degree $d=\nu(2g-2)$.  It is easy to write equations for $\phi_{\nu}(C)$---indeed, equations for the embeddings of hyperelliptic curves by more general linear systems.  Here we do so in a form convenient for our applications in the next section, following Stevens~\cite{StevensDeformations}*{pp.~137--138} and Eisenbud~\cite{EisenbudTranscanonical}. To simplify notation, we fix $\nu$ and omit it where possible. 

Let $\pi: C \rightarrow \Pro^{1}$ be the $g^{1}_{2}$ on $C$.  For convenience, write $k:= \nu(g-1)$ and $e:= g+1$. Then $\phi (C)$ lies on the scroll $S = \Pro_{\Pro^{1}}(\pi_{*}(\lbpow{\omega_C}{\nu}))$, where $\pi_{*}(\lbpow{\omega_C}{\nu}) \cong \o(k) \oplus \o(k+e)$; thus $S \cong \Pro(\o \oplus \o(-e))$. 

Let $C$ be given by the affine equation $y^{2} = f(x)$, where $f(x)$ is polynomial in $x$ of degree $2g+2$ (or $2g+1$, if the point at infinity is a branch point).  Then a basis of $H^{0}(C,\lbpow{\omega_C}{\nu})$ is given by:
\begin{equation} \label{nuKhypbasis}
B_{\nu} := \{ 1, x, x^{2}, \ldots,  x^{k} ,
y, yx, yx^{2}, \ldots, yx^{k-e}\}\,.
\end{equation}
These sections are clearly independent and we check that $(k+1)+(k-e+1) = 2 \nu(g-1) - (g+1) +2 = (2\nu -1)(g-1) = h^{0}(C,\lbpow{\omega_C}{\nu})$.  Abusing notation, we use these basis elements as variables on $\Pro^{N}$.

Equations for the scroll are also classical; modern references are \cite{HarrisFirst}*{Exercise 9.11} and  \cite{ACGHCurves}*{pp. 96--100}.  Suppose $(\nu-1)(g-1)>2$, so that $k > e$.  Then the scroll equations are given by the $2 \times 2$-minors of the {\it deleted catalecticant matrix}
\begin{equation} \label{nuKhypscrolleqns}
 M:= \left(\begin{array}{lllllllll} 
yx^{k-e} & yx^{k-e-1} & \cdots & yx & | & x^{k} & x^{k-1} & \cdots & x \\
yx^{k-e-1} & yx^{k-e-2} & \cdots & y & | & x^{k-1} & x^{k-2} & \cdots & 1
\end{array} \right) .
\end{equation}
Write $I_{S}$ for the ideal generated by the $2 \times 2$-minors of $M$.  

Next, choose a set $Q$ of quadrics encoding the following equations:
\begin{equation} \label{Quadrics}
\begin{array}{lcr}
y^{2} &  = & f(x) \\
y^{2}x & = &  xf(x) \\
& \vdots & \\
y^{2}x^{2(k-e)} & = & x^{2(k-e)} f(x).
\end{array}
\end{equation}
The particular choices of quadrics used to encode these equations won't matter once these are combined with the scroll equations.  Write $I_Q$ for the ideal generated by these quadrics. The references cited above yield, in particular:

\begin{lemma} \label{heequations} Suppose the hyperelliptic cover yielding $C$ is branched at $\infty$, so  $\deg f(x) = 2g+1$.  Then $I(\phi_{\nu}(C)) = I_{S} + I_{Q}$.  That is, the ideal of $\phi_{\nu}(C)$ is given by the scroll equations (the $2 \times 2$-minors of $M$) together with the $2(k-e)+1$ additional quadratic equations coming from $Q$. 
\end{lemma}

\begin{example} \label{wimanfourexample} To illustrate Lemma~\ref{heequations}, let's find the bicanonical equations of $\wiman{4}$ which, as we'll see in (\ref{wimanequation}) is the curve given by $y^2 = x^9 - 1$.

Here $k =6$, $e=5$, and $N = 8$.  We coordinatize $\Pro^{8}$ as follows:
\begin{displaymath}
 \begin{array}{ccccccccccccccccc}
1 & : & x & : & x^2 & : & x^3 & : & x^4 & : & x^5 & : & x^6 & : & y & : & yx\\
a & : & b & : & c   & : & d   & : & e  & : & f & : & g & : & h & : & i 
\end{array}
\end{displaymath}

The deleted catalecticant matrix is 
\begin{displaymath} M :=
\left( \begin{array}{lccccccc}
yx & | & x^6 & x^5 & x^4 & x^3 & x^2 & x \\ 
y & |  & x^5 & x^4 & x^3 & x^2 & x & 1
\end{array} \right)
=
\left( \begin{array}{cccccccc}
i & | & g & f & e & d & c & b \\ 
h & | & f & e & d & c & b & a
\end{array} \right),
\end{displaymath}
yielding (see \codesamplelink{scrollequations}{6} for the relevant \mtwonocite\ commands)
\begin{displaymath}
\begin{array}{c}
 I_{S} = (-gh+fi,-fh+ei,-f^2+eg,-eh+di,-ef+dg,-e^2+df,-dh+ci, \\
-df+cg,-de+cf,-d^2+ce,-ch+bi,-cf+bg,-ce+bf,-cd+be,-c^2+bd, \\
-bh+ai,-bf+ag,-be+af,-bd+ae,-bc+ad,-b^2+ac)
\end{array}
\end{displaymath}

We may encode the equation $y^2 = x^9 - 1$ as $h^2-dg+a^2$, the equation $y^2 x = x^{10} - x$ as $hi - eg + ab$, and the equation $y^2 x^2 = x^{11} - x^2$ as $i^2-fg+ac$.  Then the ideal we seek is 
\begin{equation} \label{wimanfourideal}
\begin{array}{c}
(-gh+fi,-fh+ei,-f^2+eg,-eh+di,-ef+dg,-e^2+df,-dh+ci,-df+cg,-de+cf, \\
-d^2+ce,-ch+bi,-cf+bg,-ce+bf,-cd+be,-c^2+bd,-bh+ai,-bf+ag,-be+af, \\
-bd+ae,-bc+ad,-b^2+ac,h^2-dg+a^2,hi-eg+ab,i^2-fg+ac).
\end{array}
\end{equation}
\end{example}

\subsection*{Bicanonical multiplicity-freeness of Wiman curves}
We write $\wiman{g}$ for the \textit{Wiman curve of type I in genus $g$}. These curves are named for Anders Wiman who, in 1895 in his first published paper\footnote{His last appeared 59 years later!}~\cite{Wiman}, showed that $\wiman{g}$ has the cyclic automorphism group of largest order $4g+2$ amongst all smooth curves of genus $g$.

The curve $\wiman{g}$ is the smooth hyperelliptic curve given by the affine equation
\begin{equation}\label{wimanequation}
	 y^{2} = x^{2g+1} -1\,.  
\end{equation}
It is often convenient to think of $\wiman{g}$ as the hypersurface in the weighted projective space $\Pro(1,g+1,1)$ given by $y^{2} = x^{2g+1}z -z^{2g+2}$ and, when we do, we call $[1:0:0]$ the branch point at infinity.  As we have already remarked,  $\Aut(\wiman{g})$ is cyclic of order $4g+2$: fixing a primitive $\thst{(4g+2)}{nd}$ root of unity $\zeta$ determines a generator $\sigma \in \Aut(\wiman{g})$ that acts with weight $2g+1$ on $y$ (that is, as $-1$) and with weight $2$ on $x$.  The key facts we will need about this action are summarized in:
\begin{proposition} Let $B_2$ be the basis of $V=H^{0}(\wiman{g},\lbpow{\omega}{2})$ given by (\ref{nuKhypbasis}) and let $T$ be the corresponding torus in $\SL(V)$. \begin{enumerate}
	\item $\Aut(\wiman{g})$ fixes the branch point at infinity.
	\item The elements of $B_2$ are eigenvectors for the action of $\sigma$ on $H^{0}(\wiman{g},\lbpow{\omega}{2})$ with distinct powers of $\zeta$ as eigenvalues.
\end{enumerate}
Hence, the bicanonical model of $W_g$ is multiplicity free and the torus $T$ determines stability for it.
\end{proposition}
\begin{proof} Since $\sigma$ clearly fixes the point at infinity, the description of its action above makes the first statement clear and shows that its action on $H^{0}(\wiman{g},\lbpow{\omega}{2})$ in the basis $B_2$ of (\ref{nuKhypbasis}) is by $\zeta^{2i}$ on $x^i$ for $i = 0, 1, \ldots, k = 2g-2$ and by $\zeta^{2i+2g+1}$ on $yx^i$ for $i = 0, 1, \ldots, k-e = g-3$. This gives the second statement and it, in turn, shows that the elements of $B_2$ span invariant lines on each of which $\Aut(\wiman{g})$ acts by a different character. From this, the final claims follow immediately.
\end{proof}

\section{Results and predictions from the log minimal model program} \label{predictionssection}

In this section, we recall some results from completed stages of the log minimal model program for the moduli stack $\Mgbar$ of stable curves and its coarse moduli scheme $\mgbar$, and some predictions about conjectural stages.  For a more detailed over view of this program, and on the moduli spaces and divisor classes discussed here, see \cites{HassettHyeonLogCanonical, HassettHyeonFlip}. Because it makes the coefficients involved slightly simpler, we'll state predictions in terms of divisors on the associated moduli stack $\Mgbar$.

The goal of the log minimal model program is to construct the birational models $\Mgbar(\alpha)$ determined by the rays $r_\alpha$ spanned by the classes $K_{\Mgbar}+ \alpha \delta$ in the Neron-Severi group of the moduli stack $\Mgbar$ of stable curves where $\delta$ is the divisor class determined by the locus of nodal curves and $1 \ge \alpha \ge 0$. In view of Mumford's formula, $K_{\Mgbar}=13\lambda -2 \delta$ in which $\lambda$ is the Hodge class, the ray $r_\alpha$ is also spanned by the class $s_\alpha\lambda - \delta$ with \emph{slope} $s_\alpha= \frac{13}{2-\alpha}$.  For $\alpha$ near $1$, the ray $r_\alpha$ is ample and $\Mgbar(\alpha) = \Mgbar$, and then, as $\alpha$ decreases, there is a discrete set of critical values of $\alpha$ or $s_\alpha$ at which the model $\Mgbar(\alpha)$ changes, until at $\alpha = 0$ we arrive as the canonical model of $\Mgbar$. Understanding the birational alterations involved and constructing the models $\Mgbar(\alpha)$ becomes more delicate as $\alpha$ decreases.

Predictions about the program can be obtained by assuming a conjectural description of the nef cone of $\Mgbar$ called the $F$-conjecture in \cite{GibneyKeelMorrison} and using consequences of it deduced there (in  particular, Theorem 2.2 and Proposition 6.1). These results lead to predictions for the successive smaller critical values $\alpha$, for the loci of stable curves contracted at each such value, and for the non-stable singularities  that replace them. Details for $\alpha = \frac{7}{10}$, for example, may be found in \cite{HassettHyeonLogCanonical}*{pp. 20--21}. For small $g$, only larger values of $\alpha$ arise and the birational alteration involved may be special---for example, loci that are usually of higher codimension  may be divisors for small $g$.  

On the other hand, a calculation due to Mumford \cite{MumfordEnseignement}*{\S 5} (see also, \cite{MorrisonGIT}*{Lemma 21}) determines the slope $s(\nu,m)$ of the polarization with which GIT naturally endows the locus of $\nu$-canonically embedded nodal curves the Hilbert scheme linearized in degree $m$ and shows that the polarization on the $\nu$-canonical the Chow scheme has slope $s(\nu,\infty)$ given by letting $m\to \infty$. His formula for the polarization is:
\begin{equation} \label{polarizationequation} 
	\left(\frac{m-1}{2 \nu -1}\right) \left(\bigl( 6 \nu^{2}m-2\nu m-2\nu+1 \bigr) \lambda
	-\frac{\nu^2 m}{2}\delta\right).
\end{equation}
When one of these GIT slopes equals a critical $s_\alpha$, the corresponding GIT quotient is a natural candidate for $\Mgbar(\alpha)$. The goal then becomes to see that the stable (and, in some cases, semistable) orbits of this quotient parameterize the expected class of curves.

Table \ref{predictiontable} summarizes the first few stages of the program. The $g_{\text{min}}$ column gives the minimal genus below which the general picture indicated must be modified in some way. Since most of these modifications will not be significant in applying our results, we pass over them. 

\def\hhh#1{\leavevmode\hbox to16pt{\hfil #1 \hfil}}
\begin{center}
	\centerline{}
			{
			\renewcommand{\arraystretch}{1.4}
			\tabcolsep=3pt
			\hyphenpenalty=10000
\begin{tabular}{|m{1.5cm}|m{.4cm}|m{.4cm}|m{.6cm}|m{3.7cm}|m{2.6cm}|m{1.9cm}|m{2.7cm}|}
			\hline 
			{\scshape status} & \hhh{$\alpha$} & $s_{\alpha}$ & $g_{\text{min}}$ & {\scshape locus contracted} &  {\scshape replacement singularities} & {\scshape alteration} & {\scshape git parameters giving alteration}\\ 
			\hline
			known & \hhh{$1$} & $13$ & --- & ---  & --- &--- &  $\nu \ge 5$, $ m\gg 0$\\
			\hline
			known & \hhh{$\frac{9}{11}$} & $11$ & 3 & elliptic tails & cusps & contraction &$\nu$\,=\,$3$\,or\,$4$, $ m\gg 0$\\
			\hline
			known & \hhh{$\frac{7}{10}$} & $10$ & 4 & elliptic bridges & tacnodes & flip &$\nu = 2$, $ m\gg 0$\\
			\hline
			predicted & \hhh{$\frac{2}{3}$} & $\frac{39}{4}$ & 5 & {Weierstrass genus $2$ tail} & {ramphoid cusps} & {flip} &{$\nu = 2$, {$ m = 6-\epsilon$}}\\
			\hline
			predicted & \hhh{$\frac{19}{29}$} & $\frac{29}{3}$ & 5 & {general genus $2$ tail} & {ramphoid cusps} & {contraction} &{$\nu = 2$, {$ m = 4.5-\epsilon$}}\\
			\hline
			\end{tabular}			
			}
			
			\centerline{}
			\centerline{\scshape Table \ref{predictiontable}~ Known and predicted stages of the log minimal model program}
	\refstepcounter{equation}\label{predictiontable}
\end{center}

Interpreting the implications of each row in the table for the GIT problems listed in the last column requires a bit of explanation. Fix a row of the table and consider the closure of the $\nu$-canonical locus in the $\nu$-canonical Hilbert scheme, linearized in degree $m$, for the values given in that row. Then, curves in the sublocus contracted in that row (and in subloci contracted in any higher rows) should have unstable Hilbert points. And, \emph{ceteris paribus}, curves with the replacement singularities in that row (and those in any higher rows, and nodes) but no others should have stable Hilbert points. The italicized proviso in this last statement is that such curves must not exhibit some other destabilizing geometric feature. A careful description of what these features are is necessary in constructing the quotient but not in dealing with the pointwise stability calculations in our examples, so we may, and will, omit giving one.

Each row also comes with an implicit limiting case obtained by sending $m \to \infty$ in the known rows and by sending $\epsilon \to 0$ in the predicted rows. In the former case, the limit is the corresponding $\nu$-canonical Chow quotient and in the latter it is the $\nu$-canonical Hilbert quotient linearized in the limiting degree $m$. In this limiting GIT problem, we expect both curves in the locus contracted in the fixed row and curves with the replacement singularities in the fixed row to be strictly semi-stable---again, absent any other destabilizing geometric feature. In particular, in the predicted rows, we expect the limiting value of $m$ to give a linearization lying on a VGIT wall in the sense of \cites{ThaddeusVGIT, DolgachevHu} and the contracted and singular curves to change from stable to unstable and from unstable to stable, respectively, as the degree $m$ descends across this value. 

\begin{remark}\label{genusthreeremark} We note one additional low genus result. Hyeon and Lee \cite{HyeonLeeGenusThree}*{Proposition 19} have shown that in genus $3$, when $\alpha < \frac{17}{28}$, the locus of hyperelliptic curves (which is a divisor in genus $3$ only) is contracted in $\Mgbar(\alpha)$. The corresponding slope $s_\alpha < \frac{28}{3}$ arises by taking $\nu=2$ and $m=\frac{9}{4}$ in \ref{polarizationequation}. Our approach to interpreting stability with respect to a fractional linearization of this type is outlined in Remark~\ref{mumremark}. Thus, here we expect hyperelliptic curves to be $m$-stable for $m>\frac{9}{4}$, $m$-strictly semistable for  $m=\frac{9}{4}$, and $m$-unstable for $m<\frac{9}{4}$.   
\end{remark}

\section{Results} \label{resultssection}
In this section, we summarize our computational results. Both the source code files used and the resulting output can be found at  \neturltilde{http://www.math.uga.edu/~davids/gs/gs.html}{http://www.math.uga.edu/$\sim$davids/gs/gs.html}.  

We prove that the bicanonical genus $3$ Wiman curve $\wiman{3}$ is unstable for $m=2$ and stable for $m \geq 3$ (which matches the predictions based on \cite{HyeonLeeGenusThree}).  We prove that the bicanonical genus $4$ Wiman curve $\wiman{4}$ is stable for $m \geq 2$, and that the genus $5$, $6$, $7$, and $8$ Wiman curves are stable for small values of $m$, which matches predictions of Hassett and Hyeon.  We study a specific genus $5$ curve with an elliptic bridge; this is unstable for all finite $m$, but Chow strictly semistable.  We also prove that a specific genus $5$ nodal curve with a genus two tail attached at a Weierstrass point is unstable for $m<6$, strictly semistable for $m=6$, and stable for $m=7$.  For $m=6$ it specializes to a curve with a ramphoid cusp, which is also strictly semistable.  We also study a specific genus $5$ nodal curve with a genus two tail attached at a non-Weierstrass point, and find that it is unstable for $m < 4.5$ and stable for $m=5,6,7$.  Finally, we study a specific genus $4$ ribbon; this is unstable for all finite $m$, but Chow strictly semistable.

\subsection*{Smooth Wiman curves}

\begin{example}[The genus 3 Wiman curve.] \label{genusthreewiman} Genus 3 bicanonical curves are not explicitly covered by Section~\ref{wimansection}, since $(\nu-1)(g-1)=2$, or equivalently, $k=e$.  But we can stretch the algorithm there to cover this case, too: instead of the curve lying on a scroll given by a deleted catalecticant matrix, in genus 3, the curve lies on a cone over the rational normal curve given by a catalecticant matrix.  To this we can add a quadric encoding $y^2 = f(x)$, yielding: the ideal for $\wiman{3}$ in $K[a,b,c,d,e,f]$ is $(ac-b^2,ad-bc,ae-bd,bd-c^2,be-cd,ce-d^2,f^2-ab+e^2)$.

In this case, it is not necessary to use Monte Carlo methods. When we compute the full state polytope, \gfan\ finds 4615 initial ideals, taking about $1.4$mb to describe. Interestingly, while $I$ is generated by quadrics, some of the initial ideals have much higher regularity---one of the initial ideals has a generator of degree 19.  (The Gotzmann number $m_{P}$ for the Hilbert polynomial $8t-2$ is 26.)  We find that $\wiman{3}$ is unstable for $m=2$, and stable for $m\geq 3$.  For $m=3, \ldots, 12$ we computed $\State_{m}(I)$; for $m=13$ through $27$ we found random sets of 6 initial ideals whose $m^{th}$ Hilbert points contained the barycenter in their convex hull.  By the Parabola Trick, this implies $\wiman{3}$ is stable for all $m\ge 3$. This corroborates the predictions of Hyeon and Lee mentioned at the end of Section~\ref{predictionssection}.  

We also computed the worst \oneps when $m=2$.  The proximum is $p = (\frac{12}{5},\frac{12}{5},\frac{12}{5},\frac{12}{5},\frac{12}{5},\frac{10}{5})$.  We first computed this using the Maple package \texttt{Convex} (\cite{Convex}); see \codesamplelink{mapleprox}{7}.  Later, we learned that version 1.1.2 of \texttt{Convex} had a bug in its proximum function, and checked its answer by verifying the Karush-Kuhn-Tucker conditions (cf. \neturltilde{http://www.math.uga.edu/~davids/gs/wiman3/kkt.pdf}{http://www.math.uga.edu/{$\sim$}davids/gs/wiman3/kkt.pdf}.)  Then $p - \bary_{2} = (\frac{1}{15},\ldots,\frac{1}{15},-\frac{1}{3})$, and therefore the worst \oneps is one that scales the span of the rational normal curve with equal weights and scales the cone with complementary weight. 

Using the $\MUm$ function from \cite{HassettHyeonLee} and $w = (10,10,10,10,10,12)$, we get $\MUm(I,w,2)=-4$ and $\MUm(I,w,3)=24$.  Interpolating using Lemma~\ref{hhlemma}, this gives $\mu([\wiman{3}]_m,\lambda) = 4(m-1)(4m-9)$.  Thus, with respect to this \oneps, $[C]$ is $m$-stable, strictly semi-stable and unstable for $m$ respectively greater than, equal to and less than~$\frac{9}{4}$, confirming the predictions discussed in Remark~\ref{genusthreeremark}.

Finally, we wish to observe that in this example, the state polytope admits a Minkowski sum decomposition: by \cite{SturmfelsGrobner}*{p. 16. Ex. 4}, we have $\State_{2}(I) = \State_{2}(I_S) + \State_2(I_Q)$.  We do not know of a similar decomposition for our other examples, or have reason to expect one.  
\end{example}

\begin{example}{\scshape The genus $4$ Wiman curve $\wiman{4}$}.\label{wimanfourcalc}  We computed the ideal of $\wiman{4}$  in (\ref{wimanfourideal}).  The Hilbert polynomial $12t-3$ is and its Gotzmann number, using Lemma \ref{gotzmannformula}, is $63$.

We have not been able to compute the full Gr\"{o}bner fan of this ideal.  To get a sense of how large this might be, we computed state polytopes for related embeddings of lower degree.  On a genus 4 hyperelliptic curve, $K \sim 6P$, where $P$ is any Weierstrass point, so it is natural to compute state polytopes for linear systems of the form $| K + nP | $.

\vskip12pt	
\begin{center}
	\begin{tabular}{|c|r|r|r|}
		\hline
		Linear system & \# of initial ideals & CPU time & \gfannocite\ output file size \\
		\hline
		$|K+3P|$ & $3{,}862$ & $3$ minutes & $1$~\hbox to 16pt{mb\hfil}\\
		\hline
		$|K+4P|$ & $283{,}221$ &  $370$ minutes  & $101$~\hbox to 16pt{mb\hfil}\\
		\hline
		$|K+5P|$ &  $20{,}694{,}486$ & $30$ days & $25$~\hbox to 16pt{gb\hfil}\\
		\hline
	\end{tabular}
\vskip9pt	
\refstepcounter{equation}\label{wimanfourtable} {\scshape Table~\ref{wimanfourtable}. Growth of complexity of state polytopes of $\wiman{4}$}
\end{center}

Table~\ref{wimanfourtable} makes it clear that computing the full Gr\"{o}bner fan for the linear system $|2K| =| K+6P|$ is out of reach.  (Note: $|K|$ is not an embedding, since $C$ is hyperelliptic.  The linear systems $|K+P|$ and $|K+2P|$ are also not embeddings.)

Next we turned to our Monte Carlo strategy.  Using random weights, we were able to establish that $\wiman{4}$ is stable for $m$ from $2$ to $7$.   Next, we examined the output from the $m=2$ calculation more closely.  Here $336$ random weights were used before stability was established.  We chose the $25$ initial ideals whose second Hilbert points were closest to the barycenter, and randomly whittled down this set to find a set of $9$ initial ideals whose $m^{th}$ Hilbert points establish stability for $4 \leq m \leq 64$ \codesamplelink{Nspanningideals}{4}.  By the Parabola Trick, this implies stability for all $m \geq 2$, and hence Chow stability.  It seems plausible that there may be a set (perhaps even many sets) of $9$ initial ideals that establish stability for all $m \geq 2$. 

For $g \geq 4$ Hyeon predicts (private communication) that divisors of slope $\leq 9$ contract the hyperelliptic locus.  Solving $\frac{20m-3}{2m} = 9$ yields $m=\frac{3}{2}$, so the prediction is that this curve should be stable for all $m \geq 2$, which matches our findings.  Moreover, we can exhibit a \oneps that flips at $m=\frac{3}{2}$.  Let $w=(-2,-2,-2,-2,-2,-2,-2,7,7)$.  Then $\MUm(I,w,2) = 108$, $\MUm(I,w,3) = 648$, $\MUm(I,w,4) = 1620$, which gives the polynomial $\mu([\wiman{4}]_m,\lambda) = 108(m-1)(2m-3)$.  By comparing $w$ to the setup used in Example  \ref{wimanfourexample}, we see that this \oneps puts all negative weights on multiples of $x$ and all positive weights on multiples of $y$.  
\end{example}

\begin{example}{\scshape Higher genus Wiman curves}.  We applied our Monte Carlo methods to the Wiman curves $\wiman{g}$ for $g=5,6,7,8$.  The ideals of these curves may be obtained using the methods of Section~\ref{wimansection}.  Table~\ref{wimanfivetable} below shows the number of random weights required to establish stability in different degrees $m$.  We have also included the Hilbert polynomial for each example and the corresponding Gotzmann number.   Our \mtwonocite\ code tested stability by adding $4(N+1)$ initial ideals at a time, and we record this, too.  Lastly, as $m$ grew large, \texttt{polymake}, run with its default settings, returned error messages in each example; the table below includes all our error-free runs, with omitted cells indicating that we encountered a \texttt{polymake} error.  Although we are confident that such lacunae could be filled by selecting different options within \texttt{polymake}, or by using a different convex geometry software package, we did not think the game worth the candle.

\vskip12pt
	\begin{center}
	\begin{tabular}{|c|c|c|c|c|}
		\hline
 & $g=5$ & $g=6$ & $g=7$ & $g=8$ \\ \hline
$P(t)$ & $16t-4$ & $20t-5$ & $24t-6$ & $28t-7$\\ \hline
$m_{P}$ & $116$ & $185$ & $270$ & $371$\\ \hline
Initial ideals per round & $48$ & $60$ & $72$ & $84$\\ \hline\hline
& \multicolumn{4}{c}{\hfil Ideals needed to establish stability}\vline \\\hline
		$m=2$ &   $336$ & $480$ & $1224$ & $1344$\\
		\hline
		$m=3$ &  $144$ & $300$& $360$ & $420$\\
		\hline
$m=4$ &  $48$ & $180$& $144$ & $168$\\
		\hline
$m=5$ &  $96$ & $120$ &  $144$ & $252$ \\
		\hline 
$m=6$ & $96$ & $60$ & $216$ & $168$ \\
		\hline
$m=7$ &  $48$ & $120$ & $72$ & $168$ \\
	\hline
$m=8$ &  $48$ & $60$ & \multicolumn{2}{c}{} \\
\cline{1-3}
$m=9$ &  $96$ &  \multicolumn{3}{c}{}  \\
\cline{1-2}
\end{tabular}
\vskip9pt	
\refstepcounter{equation}\label{wimanfivetable} {\scshape Table~\ref{wimanfivetable}. Number of random weights required to establish stability.}
\end{center}

	Table~\ref{wimanfivetable} suggests that it is easier to establish stability for larger values of $m$ than for $m=2$. For $g=5$, we also checked some slightly larger values of $m$.  Rather than generating new sets of random initial ideals, we checked that the initial ideals used in degree $7$ establish stability for degrees $4$ through $11$.  Above degree $11$, however, \texttt{polymake} returned error messages once again.    

	We expect that the whittling procedure used in Example~\ref{wimanfourcalc} could also be applied to prove stability in all degrees for the genus $5$ curve, but we have not made a systematic attempt to do so. 
\end{example}

\subsection*{Reducible curves}
\begin{example}{\scshape A genus $5$ curve with an elliptic bridge}. \label{ellbr}
An \textit{elliptic bridge} is a genus 1 subcurve that meets the rest of the curve at two nodes.  Hassett and Hyeon show that this is unstable for finite $m$, but Chow strictly semistable.  We compute an example illustrating their findings.

We build an example as follows: let $\wiman{2}$ be the Wiman curve of genus 2, $P$ the point at infinity.  Let ${E}$ be the elliptic curve given by the equation $y^2 = x^3 -x$ (Kulkarni calls curves like this the Wiman curves of type II), and let $Q = [1:0:0]$, $R = [0:0:1]$.  Then ${C} = \wiman{2} \cup_{P=Q} {E} \cup_{R=P} \wiman{2}$.  We chose coordinates for the bicanonical embedding of ${C}$ as follows:
\begin{displaymath}
\begin{array}{ccccccccccccccccccccccc}
yx & : & xz^3 & : & x^2z^2 & : & x^3z & : & x^4 &  &    &   &      &   & x^4 & : & x^3z & : & x^2z^2 & : & xz^3 & : & yx \\
   &     &    &   &      &    &        &  & x^3z & : & yxz & : &x^2z^2  & : & xz^3 &    &       &   &       &   &    &   &      \\
a & :   & b & : &  c  &  : & d     & : & e   &  : & f  & : & g  & : &  h    & : & i     & : & j      & : & k & : & l
\end{array}
\end{displaymath}
We know that $\omega_{{C}}^{2} |_{\wiman{2}} \cong \omega_{\wiman{2}}^{2}(2P)$, and $\omega_{{C}}^{2} |_{{E}} \cong \omega_{{E}}^{2}(2Q+2R) \cong \mathcal{O}_{{E}}(2Q+2R)$.  Equations for such hyperelliptic curves and  linear systems follow from the results in \cite{StevensDeformations} and \cite{EisenbudTranscanonical} but we obtained them using \Magmanocite \cites{BosmaCannonPlayoust, Magma}: see \codesamplelink{curveequations}{5} and \codesamplelink{ellipticbridgeideal}{1.1}.  

Consider the genus $2$ subcurve embedded in the span of $h,i,j,k$ and $l$.  The point at infinity maps to $[1:0:0:0:0]$.  The ideal of $\wiman{2}$ in $\Pro^4$ is given by 
\[ (l^2-hi+k^2, i^2-hj, ij-hk, j^2-ik),
\]
and the generator of $\Aut(\wiman{2})$ acts on $h$--$l$ with weights $8,6,4,2$ and $7$.  Then the ideal of $\wiman{2}$ in $\Pro^{11}$ is obtained by adding $(a,b,c,d,e,f,g)$ to this ideal, and the $\Aut(\wiman{2})$-action is extended to $\Span \{a,b,c,d,e,f,g\}$ by giving these weight $8$. 

We can get equations and automorphisms for the other tail by symmetry.  The equations and automorphisms of the elliptic curve component are also easy to find.  We obtain:
	\begin{align*}%
{I}= \quad & (a^2-de+b^2, bd-e^2, be-cd, d^2-ce, f,g,h,i,j,k,l)\\
\cap~ &(g^2-eh,f^2-eg+gh,a,b,c,d,i,j,k,l)\\
 \cap~ &(l^2-hi+k^2, i^2-hj, ij-hk, j^2-ik,a,b,c,d,e,f,g)\,.
	\end{align*}%

The following three matrices in $\GL(12)$ fix ${I}$:
{\renewcommand{\arraystretch}{1.3}
\begin{displaymath}%
\begin{array}{c}
{D}_{1} = D(\zeta_{10}^{7},\zeta_{10}^{2},\zeta_{10}^{4},\zeta_{10}^{6},\zeta_{10}^{8},\zeta_{10}^{8},\zeta_{10}^{8},\zeta_{10}^{8},\zeta_{10}^{8},\zeta_{10}^{8},\zeta_{10}^{8},\zeta_{10}^{8}) \\
{D}_{2} = D(\zeta_{10}^{8},\zeta_{10}^{8},\zeta_{10}^{8},\zeta_{10}^{8},\zeta_{10}^{8},\zeta_{10}^{8},\zeta_{10}^{8},\zeta_{10}^{8},\zeta_{10}^{6},\zeta_{10}^{4},\zeta_{10}^{2},\zeta_{10}^{7}) \\
{D}_{3} = D(-1,-1,-1,-1,-1,\zeta_{4},1,-1,-1,-1,-1,-1) \\
\end{array}%
\end{displaymath}%
}%
The representation of the subgroup of automorphisms generated by the three diagonal matrices above is not multiplicity free.  (This can be checked using \Magma\ or \GAP; for an example of such a calculation, see \codesamplelink{wimanweiertailmultfree}{8}.)  However, there is a fourth automorphism that flips the two ends of the chain.  It is given by the matrix ${A}$ that swaps $(a,l)$, $(b,k)$, $(c,j)$, $(d,i)$ and $(e,h)$, fixes $g$, and scales $f$ by~$\zeta_{4}$.

Let ${D}_{4} = {D}_{1}{D}_{2}{D}_{3}$, and let ${G} = \langle {D}_{4}, {A} \rangle$ be the subgroup of $\GL(12)$ generated by  ${D}_{4}$ and ${A}$.  ${G}$ is abelian of order $40$.  However, the basis diagonalizing ${G}$ is not the basis $a,\ldots,l$ above; thus, in the language of Definition \ref{multfreeX}, the torus ${T}$ scaling the variables $a,\ldots,l$ does not determine stability.
 
The basis diagonalizing ${G}$ is $a+l$, $a-l$, $b+k$, $b-k$, $c+j$, $c-j$, $d+i$, $d-i$, $e+h$, $e-h$,$f$,$g$.  We relabel these as variables $A,B,C,\ldots,L$, write ${I'}$ for the ideal in these coordinates \codesamplelink{ellipticbridgeideal}{1.1}, and write ${T}'$ for the torus scaling these variables.  

We normalize the matrices ${D}_{4}$ and ${A}$ to have determinant $1$, and change them to this basis, yielding:
{\renewcommand{\arraystretch}{1.3}
\begin{displaymath}%
\begin{array}{c}
{D}_{4}' = D(\zeta_{240}^{3},\zeta_{240}^{3},\zeta_{240}^{123},\zeta_{240}^{123},\zeta_{240}^{171},\zeta_{240}^{171},\zeta_{240}^{219},\zeta_{240}^{219},\zeta_{240}^{27},\zeta_{240}^{27},\zeta_{240}^{207},\zeta_{240}^{147})\\
{A}'  = D(\zeta_{240}^{5},\zeta_{240}^{125},\zeta_{240}^{5},\zeta_{240}^{125},\zeta_{240}^{5},\zeta_{240}^{125},\zeta_{240}^{5},\zeta_{240}^{125},\zeta_{240}^{5},\zeta_{240}^{125},\zeta_{240}^{65},\zeta_{240}^{5})
\end{array}
\end{displaymath}
}%
The subgroup ${G}' = \langle {D}_{4}', {A}' \rangle$ of $\SL(12)$ is multiplicity free.  

The generators of ${I}'$ are far more complicated than the generators of ${I}$, and the state polytope $\State_{{T}'}({I}')$ is likely to be correspondingly more complicated than $\State_{{T}}({I})$. Since we found that ${I}$ already has $500{,}094$ initial ideals, generated in degrees $2$ through $9$, we decided to work in its coordinate system as much as possible.  

The elliptic bridge is known to be bicanonically unstable in all degrees $m \geq 2$ by \cite{HassettHyeonFlip}*{\S 10} where the \oneps with weights $(2,2,2,2, 2,1,0,2, 2,2,2,2)$ is shown to be destabilizing.  As a check, we computed $\MUm({I},w,2) = -12$, $\MUm({I},w,2) = -24$, and $\MUm({I},w,4) = -36$; hence $ \mu([{C}]_m,\lambda) = -12(m-1)$. As a further check, we verified in the extreme degrees, $2$ and $9$, that the corresponding state polytope did not contain the barycenter.

The coefficient of $m^2$ in $ \mu([{C}]_m,\lambda) $ is $0$. Thus, this \oneps is not Chow destabilizing. Indeed, this elliptic bridge is Chow strictly semistable, as \cite{HassettHyeonFlip}*{\S 11} shows by deformation theoretic arguments.  To check this, we needed to switch to the torus ${T}'$, since this torus determines stability. We did not attempt to compute the full Chow polytope of $I'$, but instead performed a Monte Carlo calculation which found ideals that confirmed strict Chow semistability.

\end{example}

\begin{example}{\scshape A genus $5$ curve with a genus $2$ tail attached at a Weierstrass point}.  \label{g2wp} Here we consider an example of a nodal genus $5$ curve that has a genus $3$ component and a genus $2$ component (hence a genus $2$ tail) where the node is a Weierstrass point of the genus $2$ component.  Hassett and Hyeon predict that such a curve is stable for $m>6$, semistable for $m=6$, and unstable for $m<6$.  Our calculations confirm this prediction.  

We build our $C$ by letting $\wiman{3}$ be the Wiman curve of genus $3$ with $P$ its the branch point at infinity, letting $\wiman{2}$ be the genus $2$ Wiman curve with  $Q$ its branch point at infinity, and setting $C = \wiman{3} \cup_{P=Q} \wiman{2}$.   The linear series $\omega_{C}^{2}$ is very ample, and the image of $C$ under the corresponding morphism $\phi$ is a degree $16$ curve in $\Pro^{11}$. 
We know that $\restrictedto{\omega_{C}^{2}}{\wiman{3}} = \omega_{\wiman{3}}^{2}(2P)$, and $\restrictedto{\omega_{C}^{2}}{\wiman{2}} = \omega_{\wiman{2}}^{2}(2Q)$.   Once again we used \Magmanocite to obtain the corresponding equations; see \codesamplelink{curveequations}{5}. 

For $\wiman{3}$, coordinatize $\Pro^7$ using the variables $a$--$h$, and map $\wiman{3} \rightarrow \Pro^{7}$ by
\begin{displaymath}
\begin{array}{ccccccccccccccc}
yx^6&   & yzx^5  &   & z^5x^5 &   & z^4x^6&   &z^3x^7  &   &z^2x^8 &  & zx^9 &  &  x^{10}    \\
a & : & b  & : &  c  & : & d  & : & e  & : & f & : & g  & : & h
\end{array}
\end{displaymath}
Then $P$ maps to $[0:0:0:0:0:0:0:1]$.  The ideal of $\wiman{3}$ in $\Pro^{7}$ is given by
\begin{displaymath}
\begin{array}{l} (-g^2 + fh, -fg + eh, -f^2 + eg, -f^2 + dh, -ef + dg, -ef + ch, -e^2 + df, -e^2 + cg, \\ 
-de + cf, -d^2 + ce, ag - bh, af - bg, ae - bf, ad - be, ac - bd,  \\
b^2 + c^2 - fg, ab + cd - g^2, a^2 + d^2 - gh)
\end{array}
\end{displaymath}
and the generator of $\Aut(\wiman{3})$ acts on $a$--$h$ with weights $5,3,10,12,0,2,4,6$.  Then the ideal of $\wiman{3}$ in $\Pro^{11}$ is obtained by adding $(i,j,k,l)$ to the ideal above, and the $\Aut(\wiman{3})$-action is extended to $\Span \{i,j,k,l\}$ by giving these weight 6.
  
Equations and automorphisms of the $\wiman{2}$ component can be found in Example \ref{ellbr}.  We intersect the ideals of the two components to get the ideal of $C$ \codesamplelink{g2wpideal}{1.2}.  Its Hilbert polynomial is $P(m) = 16m-4$.

We check that the $\Aut(C)$-action is multiplicity-free. 
The $\Aut(C)$ representation is generated by the diagonal matrices 
\[ D(\zeta_{14}^{5},\zeta_{14}^{3},\zeta_{14}^{10},\zeta_{14}^{12},\zeta_{14}^{0},\zeta_{14}^{2},\zeta_{14}^{4},\zeta_{14}^{6},\zeta_{14}^{6},\zeta_{14}^{6},\zeta_{14}^{6},\zeta_{14}^{6})
\]
and 
\[ D(\zeta_{10}^{8},\zeta_{10}^{8},\zeta_{10}^{8},\zeta_{10}^{8},\zeta_{10}^{8},\zeta_{10}^{8},\zeta_{10}^{8},\zeta_{10}^{8},\zeta_{10}^{6},\zeta_{10}^{4},\zeta_{10}^{2},\zeta_{10}^{7}).
\]  
in $\GL(12)$.  These can be normalized to elements of $\SL(12)$: 
\[ D(\zeta_{42}^{16},\zeta_{42}^{10},\zeta_{42}^{31},\zeta_{42}^{37},\zeta_{42}^{1},\zeta_{42}^{7},\zeta_{42}^{13},\zeta_{42}^{19},\zeta_{42}^{19},\zeta_{42}^{19},\zeta_{42}^{19},\zeta_{42}^{19})
\]
and 
\[ D(\zeta_{120}^{103},\zeta_{120}^{103},\zeta_{120}^{103},\zeta_{120}^{103},\zeta_{120}^{103},\zeta_{120}^{103},\zeta_{120}^{103},\zeta_{120}^{103},\zeta_{120}^{79},\zeta_{120}^{55},\zeta_{120}^{31},\zeta_{120}^{91}).
\]
The product of these two matrices is diagonal and has distinct eigenvalues; hence, the representation of the cyclic group it generates is multiplicity-free.  

We ran our Monte Carlo program for small values of $m$.  We found, as expected, that for $m=6$, $C$ is strictly semistable, and stable for $m = 7$.  As a complement, we prove deductively that $C$ is unstable for $m<6$ and strictly semi-stable for $m=6$.  Let $\lambda$ be the \oneps that acts with weights $w=(6,6,6,6,6,6,6,6,4,2,0,5)$. In the notation of \cite{MorrisonGIT}*{Proposition 3}, $\lambda$ has average weight $\alpha = \frac{59}{12}$ and induces a weight filtration on $\o(m)$ that looks like
\[ \underbrace{6m, 6m, ...6m}_{10m-2}, 6m-1, 6m-2, 6m-3, ..., 6,5,4,2,0,
\]
giving $w(m) = 78m^2-15m-4$.  Putting this all together gives
\[ \mu(C, \lambda)(m) = -\bigl(w(m) - mP(m)\alpha\bigr) = (\frac{2}{3})(m-1)(m-6)
\] 
with the desired roots and sign for asymptotic stability with a flip at
$m=6$. A check is provided by running the \mtwonocite\ function \texttt{MUm} on this $C$ and $\lambda$. We find that $\MUm(I,w,2)=  -32$ and $\MUm(I,w,3) = -48$ which implies by Lemma~\ref{hhlemma} that $ \mu([C],\lambda)(m) = 8(m-1)(m-6)$. The factor $12$ by which this results differs comes from the fact that \texttt{MUm} shifts the weights by $-\alpha = \frac{59}{12}$ to normalize their sum to be $0$, making it necessary to scale them by $12$ to make them integral. 
\end{example}

\begin{example}{\scshape A reducible genus $5$ curve with a ramphoid cusp}.
The calculation performed for the example above 
also establishes semistability of a curve with a ramphoid cusp as well and allows us to confirm the deformation theoretic analysis of \cite{HassettHyeonFlip}*{\S 10} at the level of ideals.

Let $C$ be the curve of the previous example, and once again let $\lambda$ be the \oneps acting with weights $w=(6,6,6,6,6,6,6,6,4,2,0,5)$.  Let $C'$ be the ideal obtained as the $\lambda$ limit of $C$ (but do not break ties, so that $C'$ is not defined by a monomial ideal).  The matrices given above for $C$ also fix $I(C')$, so $C'$ is again multiplicity free.  As shown in \cite{SturmfelsGrobner}*{Lemma 2.6}, $\State_{C'}(I)$ appears as a facet of $\State_{C}(I)$.  The barycenter is on this facet, so we also get Monte Carlo $m=6$ semistability for $C'$.

Finally, we examine $C'$ more closely.  On all but one of the generators for $I(C)$, the \oneps $\lambda$ acts with equal weights on each term, and so these are unchanged in the $\lambda$ specialization.  In particular, the $\wiman{3}$ component is unchanged in the $\lambda$ specialization.   In contrast, the generator $hi-k^2-l^2$ specializes to $hi - l^2$ under $\lambda$.  Recalling our embedding of $\wiman{2}$, we see that $hi-l^2$ encodes $x^{7}z - y^{2} x^{2}$, or $y^2 = x^5$, which is a ramphoid cusp.  Since a ramphoid cusp adds $2$ to the arithmetic genus, and since this component of $C'$ is the flat limit of the smooth curve $\wiman{2}$ of arithmetic genus 2, we see that $C'$ must be a reducible curve with two irreducible components meeting at a node.  One component of $C'$ is $\wiman{3}$, and the other component is a rational curve with a ramphoid cusp.  Note that the rational component has a $\mathbb{G}_{m}$ of automorphisms fixing the node and cusp that accounts for the $m=6$ strict semistability.  
\end{example}

\begin{example}{\scshape A reducible genus $5$ curve with a general genus $2$ tail}. \label{g2nonwp}
Here we consider an example of a nodal genus $5$ curve that has a genus $3$ component and a genus $2$ component (hence a genus $2$ tail) where the node is not a Weierstrass point of the genus $2$ component.  Hassett and Hyeon predict that such a curve is stable for $m>4.5$, semistable for $m=4.5$, and unstable for $m<4.5$.  Our calculations again confirm these predictions.

As in Example~\ref{g2wp}, we will use $\wiman{3}$ for the genus $3$ component, and $P$ its point at infinity.  For the genus $2$ component, we use a twisted model of $\wiman{2}$.  Namely, let $D$ be given by $y^2 = x^5+1$, and let $Q$ be the point $[0:1:1]$.  The hyperelliptic involution sends $Q$ to $[0:-1:1]$, and so $Q$ is not a Weierstrass point.  On the other hand, $Q$ is fixed under the automorphism $T: D \rightarrow D$ given by $(x,y,z) \mapsto ( \zeta_{5} x, y, z)$.   We used \Magmanocite to obtain equations of $D$ under the embedding $|2K + 2Q|$; see \codesamplelink{curveequations}{5}.  We coordinatize $\Pro^4$ using the variables $h$--$l$, and map $D \rightarrow \Pro^{4}$ is given by
\begin{displaymath}
\begin{array}{ccccccccc}
yz^4+z^7&   & xyz^3 +xz^6  &   & x^4z^3 &   & x^3z^4&   &x^2z^5   \\
h & : & i  & : &  j  & : & k  & : & l
\end{array}
\end{displaymath}
Then $Q$ maps to $[1:0:0:0:0]$.  The ideal of $\wiman{2}$ in $\Pro^4$ is given by 
\[ (k^2-jl, i^2-jk+2hl, hk-il, hj-ik),
\]
and $T$ acts on $h$--$l$ with weights $0,1,4,3,2$.  Then the ideal of $D$ in $\Pro^{11}$ is obtained by adding $(a,b,c,d,e,f,g)$ to this ideal, and the $T$-action is extended to $\Span \{a,b,c,d,e,f,g\}$ by giving these weight $0$. 

We intersect this with the ideal of the genus $3$ component to obtain the ideal of $C$ \codesamplelink{g2nonwpideal}{1.3}.  Its automorphism group is generated by 
\[ D(\zeta_{42}^{16},\zeta_{42}^{10},\zeta_{42}^{31},\zeta_{42}^{37},\zeta_{42}^{1},\zeta_{42}^{7},\zeta_{42}^{13},\zeta_{42}^{19},\zeta_{42}^{19},\zeta_{42}^{19},\zeta_{42}^{19},\zeta_{42}^{19})
\]
and
\[ D(1,1,1,1, 1,1,1,1, \zeta_{5}^{1},\zeta_{5}^{4},\zeta_{5}^{3},\zeta_{5}^{2}).
\]
The product of these matrices is diagonal with distinct entries, and hence $C$ is multiplicity free for the cyclic subgroup it generates.

We ran our Monte Carlo program for small values of $m$.  We found, as expected, that for $m=5,6,7$, $C$ is stable. 

Next, we studied the \oneps $\lambda$ with weights $w = (4,4,4,4, 4,4,4,4, 3,0,1,2)$ (that is, weights complementary to the weights of $T$).  We computed $\MUm(I,w,2) = -20$, $\MUm(I,w,3) = -24$, and $\MUm(I,w,4) = -12$; hence $ \mu([C]_m,\lambda) = 8(m-1)(m-\frac{9}{2})$.   Thus, this curve is unstable for $m < 4.5$.

Once again, we study the $\lambda$ specialization $C'$.  As before, most generators in the ideal are unchanged in the limit, with the exception of the last generator,  $i^2-jk-2hl$, which specializes to $i^2-jk$ under $\lambda$.  Recalling our embedding of $D$, we see that $i^2-jk$ encodes $(y+z^3)^2 - x^5 z$, which is again a ramphoid cusp, but without a hyperelliptic involution this time.
\end{example}

\subsection*{A nonreduced curve} \label{nonreducedex}
\begin{example}{\scshape A genus $4$ ribbon}.  Recall that a \emph{ribbon} is just a double structure on $\Pro^{1}$---that is, a scheme $C$ such that $C_{\red} \cong \Pro^{1}$ and $\scr{I}_{C/C_{\red}}^{2} = 0$.  Ribbons arise as limits of canonical curves, and thus it is natural to study their GIT stability alongside examples from the log minimal model program.
	
	We study the genus 4 example given in \cite{BayerEisenbudRibbons}*{p. 475}.  Let $I \subset k[a,b,c,d]$ be the ideal $ \langle ac-b^2,ad^2-2bcd+c^3 \rangle$.  This ideal admits a $\mathbb{G}_{m}$-action with weights $-3,-1,1,3$ on the variables $a,b,c,d$, respectively. Its  ideal $I$ has twelve initial ideals (listed in Table~\ref{nonreducedtable}), all generated in degrees $\leq 6$.

	For any finite $m$, one can easily check that the ribbon is Hilbert $m$-unstable using the $\lambda$ determined by the $\mathbb{G}_{m}$ above.  But we can use the approach described in Remark~\ref{chowremark} to handle all $m \geq 6$ by computing the vector $\xi_m$ associated to each initial ideal in degrees $6$ through $8$ and interpolating to obtain the polynomial representing $\xi_m$ for any $m \geq 6$ (cf. \cite{KapranovSturmfelsZelevinsky}*{p. 202}). These are also given in the table.
	
	We observe that the polytope these span, whose dimension we would expect to be $3$, lies in a plane.  The ``extra'' normal vector besides $(1,1,1,1)$ is $(-3,-1,1,3)$, the weight vector of the $\mathbb{G}_{m}$-action.  Thus, $\State_{m}(I)$ is also only two-dimensional.  For example, when $m=6$, we have that $\State_{m}(I)$ is contained in the plane defined by the equations $a+b+c+d = 306$ and $-3a -b+c+3d = -14$.  On the other hand, the barycenter does not satisfy the second equation, and is therefore outside the state polytope.

\begin{displaymath}
	\begin{array}{ll}
	\langle b^2d^2, ad^2, ac \rangle  & (m^2, 2m^2-m-1, 2m^2-4m+4, m^2+2m-3)  \\
	\\
	\langle c^4, ad^2, ac \rangle   & (m^2, 3m^2-6m+5, 6m-8, 2m^2-3m+3) \\
	\langle c^3, ac, a^2d^2  \rangle   &  (m^2+m-2, 3m^2-6m+5, 3m-2, 2m^2-m-1) \\
	\langle c^3, b^2c^2, ac, a^3d^2 \rangle &  (m^2+3m-8, 3m^2-8m+11, m+4, 2m^2+m-7) \\
	\langle c^3, b^2c^2, b^4c, ac, a^4d^2 \rangle & (m^2+6m-20, 3m^2-12m+27, 8, 2m^2+3m-15) \\
	\\
	\langle b^2, ad^2 \rangle  &(2m^2-2m+1, 3m-3, 3m^2-6m+5, m^2+2m-3) \\
	\langle bcd, b^2, abd^2, a^2d^3 \rangle  &( 2m^2-4, m+2, 3m^2-8m+10, m^2+4m-8) \\
	\langle bcd, bc^4, b^2, abd^2, a^2d^4 \rangle & (2m^2+2m-12, 6, 3m^2-12m+26, m^2+7m-20) \\
	\\
	\langle c^3, b^2 \rangle & (3m^2-6m+5, 3m-3,6m-7, 3m^2-6m+5) \\
	\langle c^3, b^2c^2, b^4c, b^6, ac \rangle  & (3m^2-12m+20, 15m-33, 8, 3m^2-6m+5)\\
	\langle c^6, bcd, bc^3, b^2, abd^3 \rangle  &(3m^2-7m+9, 6, 15m-37, 3m^2-11m+22) \\
	\langle c^5, bcd, bc^3, b^2 \rangle & (3m^2-6m+5, m+2, 10m-17, 3m^2-8m+10) 
	\end{array}
	\end{displaymath}
	\vskip3pt	
\centerline{	\refstepcounter{equation}\label{nonreducedtable} {\scshape Table~\ref{nonreducedtable}. Initial ideals of a genus $4$ ribbon.}}
\vskip6pt	

Finally, we compute the Chow polytope of $I$.  Examining the list above, we see that $\Chow(I) $ has four vertices $(1, 2, 2, 1)$, $(1, 3, 0, 2)$, $(3, 0, 0, 3)$ and $(2, 0, 3, 1)$, each approached by one of the $4$ groups of ideals above.  This quadrangle contains the barycenter of its ambient plane, so $I$ is Chow semistable.  
In summary, this ribbon is Hilbert unstable for all finite $m \geq 2$, but Chow strictly semistable.

\begin{remark}  GIT stability and semistability are open conditions.  Thus, in our previous examples, whenever we found that our example was (semi)stable for a given linearization, this implied that a general member of the same component of $\Hilbhat$ was also (semi)stable.  But since the unstable locus is closed, and multiplicity free examples are very special, the behavior of the particular example above does not indicate that a general ribbon is unstable for finite $m$.  
\end{remark}

\end{example}

\section{Future steps} \label{thoughtssection}

We hope to extend this work in several directions. First, we would like to understand other examples with geometry suggested by the log minimal model program. One such class is that of irreducible curves with a ramphoid cusp and genus at least the $g_{\min}$ of $5$ for such curves (cf. Table~\ref{predictiontable}). Another, suggested by recent work of Smyth~\cites{SmythGenusOne, SmythModular} is the class of ``elliptic triboroughs'', curves with a genus $1$ component meeting the rest of the curve in $3$ nodes. To date, we can neither find multiplicity free examples, nor show that such examples do not exist, in either class.

Although the examples here provide numerical evidence for conjectural stages of the log minimal model program, they are far from constructing any of the quotients that would be needed to verify these conjectures. Such constructions remain our ultimate objective.

As a first step, we want to prove the non-emptiness of the stable loci involved by showing that the Wiman curves of all genera are bicanonically stable. To do so, we must better understand the geometry of the initial ideals arising in our examples, with the aim of finding patterns that will allow us to replace our computational proofs of stability by deductive ones. We can identify, in the \gfannocite\ output for our small genus examples, initial ideals for which the one-parameter degeneration that produces them can be understood geometrically, in terms that do not depend on $g$. This allows us to write down an analogous degenerations for any $g$ and our goal is to use these to predict the exponent vectors of the corresponding monomial limits. Finally, we will need to be able construct enough such degenerations to prove that the convex hull of their monomial limits always contains the relevant barycenter. The proof of Chow semistability in Example~\ref{nonreducedex} can be viewed as toy model for this plan.

Second, it will be necessary to pass from the multigraded Hilbert schemes $\Hilbhat$ used here to the corresponding Grothendieck Hilbert schemes $\Hilb$, since it is quotients of the latter that naturally carry the polarizations needed to construct further log minimal models. Doing so would require, for example, showing that no codimension $1$ component of the complement of the $\ell$-nice locus in $\Hilb$ lies in the $m$-stable locus for the relevant degree $m$. Compared to the previous problem, this is, in some ways, much harder, since it requires dealing with curves exhibiting the menagerie of pathologies typical of the Hilbert scheme, and, in others, easier, since what must be checked is that such curves are \emph{not} stable. 

\newcounter{lastbib}

\linespread{.98}\normalfont\selectfont

\setcounter{lastbib}{\value{bib}}

\section*{Software Packages Referenced}

\end{document}